\documentclass[12pt]{article}
\usepackage{amsmath}
\usepackage{amssymb}
\usepackage{amsfonts}
\usepackage{eucal}
\usepackage[usenames]{color}
\usepackage{graphicx}
\numberwithin{equation}{section}
\newtheorem{Theorem}{Theorem}[section]
\newtheorem{Proposition}[Theorem]{Proposition}

\newtheorem{Remark}[Theorem]{Remark}

\begin{document}
\title{Study of direct and inverse first-exit  problems for drifted Brownian motion  with Poissonian resetting}

\author{Mario Abundo\thanks{Dipartimento di Matematica, Universit\`a  ``Tor Vergata'', via della Ricerca Scientifica, I-00133 Rome, Italy.
E-mail: \tt{abundo@mat.uniroma2.it}}
}
\date{}
\maketitle

\begin{abstract}
\noindent
We address some direct and inverse problems, for the first-exit time (FET) $\tau $ of a drifted Brownian motion  with Poissonian resetting  ${\cal X}(t)$  from an interval $(0,b)$ and the first-exit area (FEA) $A,$ namely the area swept out by ${\cal X}(t)$  till the time $\tau $; this type of diffusion process ${\cal X}(t)$ is characterized by the fact that a reset to the position $x_R $ can occur according to a homogeneous Poisson process with rate $r>0.$
When the initial position ${\cal X}(0)= \eta \in (0,b)$  is deterministic and fixed, the direct FET problem consists in investigating the statistical properties of the FET $\tau ,$
whilst the direct FEA problem studies the probability distribution of the FEA $A$.
The inverse FET problem regards the case when $\eta $ is randomly distributed in $(0,b)$ (while $r$ and $x_R $ are fixed); if
 $F(t)$ is a given distribution function on the time $t$ axis, the inverse FET  problem  consists in finding the density $g$ of  $\eta,$  if it exists, such that $P[\tau \le t ] = F(t), \ t >0.$
Several explicit examples of solutions to the inverse FET problem are provided.

\end{abstract}


\section{Introduction}
This paper continues the articles \cite{abundo:TPMS2024} and \cite{abundo:FPA2023}, which respectively concern inverse problems for the first-passage place and the first-passage time, and  the first-passage area of a one-dimensional diffusion process $\mathcal X(t)$  with Poissonian resetting in an interval $(0, + \infty).$ Here, we address the analogous direct and inverse first-exit time (FET) problems for  $\mathcal X(t)$ in an interval $(0,b).$  The direct FET problem was also treated in \cite{guoyan:24},  \cite{huang:24}, using a slight different approach and a language derived from  Physics. In the present paper, our aim is to obtain the quantities of interest as solutions of certain ordinary differential equations with boundary conditions, which in some cases we solve explicitly in terms of elementary functions. As concerns the direct FET problem, for  fixed $\mathcal X(0) = x \in (0,b)$ we investigate the statistical properties of the FET $\tau (x) $ of $\mathcal X(t)$ from the interval $(0,b),$
under the condition that $\mathcal X(0)=x ,$ that is,
$\tau (x) = \min \{t>0: \mathcal X(t) \notin (0,b) | \mathcal X(0) =x \};$  moreover, we study the probability distribution of the first-exit area (FEA) $A (x),$ namely the area swept out by ${\cal X}(t)$  till the time $\tau (x),$ and the maximum and minimum displacement of ${\cal X}(t)$ for  $t \in (0, \tau (x)).$ \par
As for the inverse first-exit time (IFET) problem, we suppose that the starting position $\eta = \mathcal X(0)$ is randomly distributed in $(0,b),$ and we denote by $\tau$ the FET of $\mathcal X(t)$ from the interval $(0,b);$ then, for a given distribution function $F(t)$ on the the time $t$ axis, the IFET problem consists in finding the probability density $g$ of  $\eta,$  if it exists, such that
the FET of $\mathcal X(t)$ from $(0,b)$ has distribution $F(t),$ namely
$P[\tau \le t ] = F(t), \ t >0;$ the density $g$ is called a solution to the  IFET problem for $\mathcal X(t).$ The analogous IFET problem for diffusions without resetting was studied in \cite{abundo:stapro13}. \par
Now, we briefly recall the definition of a one-dimensional diffusion process  with Poissonian resetting, $\mathcal X(t).$
\par
Let $X(t)$ be a one-dimensional temporally homogeneous diffusion process,
driven by the SDE:
\begin{equation} \label{diffusion}
dX(t)=  \mu(X(t)) dt + \sigma (X(t)) dW_t ,
\end{equation}
and starting from
an initial position $X(0)=\eta \in (0,b)$ (fixed or random),
where $W_t$ is a standard Brownian motion (BM) and the drift $\mu (\cdot)$ and diffusion coefficient $\sigma (\cdot)$ are regular enough functions, such that there exists a unique strong solution of the
SDE \eqref{diffusion} (see e.g. \cite{klebaner}).
\par\noindent
From $X(t)$ we construct a new process ${\cal X}(t),$ as follows.
We suppose that resetting events can occur according to a homogeneous Poisson process with rate $r>0.$
Until the first resetting event the process ${\cal X}(t)$ coincides with $X(t)$ and it evolves according to \eqref{diffusion} with ${\cal X}(0)=X(0)=\eta ;$ when the reset occurs,
${\cal X}(t)$ is set instantly to a position $x_R \in (0,b).$ After that, ${\cal X}(t)$ evolves again according to \eqref{diffusion}, starting afresh (independently of the past history) from $x_R,$ until the next resetting event occurs, and so on. The inter-resetting times turn out to be independent and exponentially distributed random variables with parameter $r.$
In other words, in any time interval $(t, t+ \Delta t),$ with $\Delta t \rightarrow 0 ^+, $ the process can pass from ${\cal X}(t)$  to the position $x_R$ with probability $r \Delta t  + o( \Delta t),$ or it can continue its evolution according to \eqref{diffusion} with probability $1- r \Delta t + o( \Delta t ).$
\par\noindent
The process ${\cal X}(t)$ so obtained is called diffusion process with Poissonian resetting,  or more simply diffusion  with  stochastic resetting;
it has some analogies with the process considered in \cite{dicre:03}, where the authors studied a M/M/1 queue with catastrophes and its  continuous approximation, namely a Wiener process subject to randomly occurring jumps at a given rate $\xi,$ each jump making the process instantly obtain the state $0.$
 Thus, the process considered in \cite{dicre:03} can be viewed as a Wiener process with resetting,
 in which a reset to the position $x_R =0$ is done, according to a homogeneous Poisson process with rate $r= \xi.$
\par
For any $C^2$ function $u(x),$ the infinitesimal generator of ${\cal X}(t)$ is given by (see e.g. \cite{abundo:FPA2023}):
\begin{equation} \label{generator}
{\cal L}u(x) = \frac 1 2 \sigma ^2(x) u''(x) + \mu (x) u'(x) +r (u(x_R) -u(x)) \equiv L u(x) +r (u(x_R) -u(x)) ,
\end{equation}
where $u'(x)$ and $u''(x)$ stand for the first and second derivative of $u(x).$
Here, $Lu(x)= \frac 1 2 \sigma ^2(x) u''(x) + \mu (x) u'(x)$ represents the  ``diffusion part'' of the generator, i.e. that concerning the diffusion process $X(t).$
\par
Direct and inverse problems for the FET and the FEA of diffusion processes are worthy of attention, since they
have interesting applications in several applied fields, for instance in biology in the context of diffusion models
for neuronal activity  (see e.g. \cite{lanska:89}, \cite{norisa:85} and the references contained in \cite{abundo:FPA2023}).
They are also relevant in
Mathematical Finance, in particular in credit risk modeling (see e.g.
\cite{jackson:stapro09}); other applications can be found e.g.
in queuing theory (see e.g. the discussion in \cite{abundo:stapro12});
for a review concerning Brownian motion with resetting in physics and computer science, see e.g \cite{maj07}. For more about inverse passage problems concerning jump-diffusions without resetting, see e.g.
\cite{abundo:CC22}, \cite{abundo:saa20IFPP}, \cite{abundo:saa19}, \cite{abundo:mathematics}, \cite{abundo:LNSIM},  \cite{abundo14b}, \cite{abundo:MCAP2013}, \cite{abundo13b}, \cite{abundo:stapro13}, \cite{abundo:stapro12}, \cite{abundo:pms00}, \cite{lefeb:19}, \cite{lefeb:22}, \cite{tuckwell:76}.
\par
At our knowledge,  extensions of inverse FET  problems to diffusions with resetting have not been treated  in the literature, yet. Thus, the aim of the present article is to study these types of problems for this class of
processes, by characterizing the quantities of interest as solutions of suitable ordinary differential equations with boundary conditions. \par\noindent
Of course, one can study FET problems of a diffusion  with resetting in any interval $(a,b);$  we have taken $a=0,$ only for the sake of simplicity.\par
Although we obtain explicit formulae only in the special case  when ${\cal X}(t)$ is a drifted Brownian motion with resetting, in
principle analogous calculations can be developed  for a general one-dimensional diffusion with resetting, whose underlying diffusion is driven by the SDE \eqref{diffusion}.
It suffices to replace, in all the differential equations, the generator of drifted Brownian motion, i.e. $Lu(x) = \frac 1 2 u''(x)+ \mu u'(x),$ with the diffusion part
of the generator $\mathcal L$ given by \eqref{generator}, that is, $Lu(x)=
\frac 1 2 \sigma^ 2(x) u''(x)+ \mu (x) u'(x);$ of course the computations turn out to be more complicated.
\par
The paper is
organized as follows. Section 2 treats
the direct FET problem when
${\cal X}(t)$ is drifted Brownian motion with resetting:
suitable differential problems are stated for the exit probability of ${\cal X}(t)$ from the left end of the interval $(0,b),$ for the Laplace transform and the moments of the FET $\tau(x)$ and  the FEA $A(x),$  and for the joint moment of $\tau(x)$ and $A(x);$ finally the probability distributions of the maximum and minimum displacement of
${\cal X}(t)$
till  the FET $\tau(x)$ are investigated. In Section 3 the inverse first-exit time (IFET) problem for
drifted Brownian motion with resetting is studied and several explicit examples of solutions are provided. Finally, Section 4 contains conclusions and final remarks.

\section{The direct first-exit time problem for drifted Brownian motion with resetting}
In this Section,  $\mathcal X(t)$ is drifted Brownian motion with resetting, also called Wiener process with resetting, and  we  analyze the direct first-exit time problem in $(0,b);$ therefore,
we suppose that the underlying diffusion is  $X (t)= x + \mu t + W_t, \ x \in (0,b).$
Then, for any $C^2$ function $u(x),$ the infinitesimal generator of ${\cal X}  (t)$ is  (see \eqref{generator}):
\begin{equation} \label{generator1}
{\cal L}u(x) =  Lu(x) + r (u(x_R) -u(x)) = \frac 1 2 u''(x) + \mu  u'(x) +r (u(x_R) -u(x)) ,
\end{equation}
where $Lu(x)= \frac 1 2 u''(x)+ \mu u'(x)$
(we omit to specify the dependence of $L$ and ${\cal L}$ on $\mu,$ for simplicity). \par\noindent
We will find explicit solutions, in terms of elementary functions only, to
some differential problems with boundary conditions (see e.g. \cite{abundo:TPMS2024}, \cite{abundo:FPA2023}), concerning
the exit probability of $\mathcal X (t)$ from the left end of the interval $(0,b),$
and
the statistical properties of the FET $\tau _ \mu (x)$ and the FEA $A _ \mu (x),$ where
\begin{equation}
\tau _ \mu (x)= \min \{t>0: \mathcal X(t) \notin (0,b) | \mathcal X(0) =x  \}, \ x \in (0,b),
\end{equation}
denotes the first-exit time (FET) of $\mathcal X  (t)$ from the interval $(0,b),$ under the condition that $\mathcal X(0)=x \in (0,b),$
and
\begin{equation}
A _ \mu (x) = \int _0 ^ {\tau (x)} \mathcal X(t) dt
\end{equation}
is the first-exit area (FEA), i.e. the area swept out by ${\cal X}  (t)$  till the time $\tau _ \mu (x) ,$ under the condition that $\mathcal X(0)=x \in (0,b).$
Moreover, we will study the probability distributions of the maximum and minimum displacement of $\mathcal X(t)$ till the FET.
For simplicity of notations, we omit to explicit the dependence of all the quantities on $r$ and $x_R.$ \par\noindent
Since the moments of the FET $\tau _ \mu (x)$ are finite (this will be shown in Subsection 2.3), then the moments of $A _ \mu (x)$ are also finite, because $A _ \mu (x) \le b \tau _ \mu (x),$ and so
$E[(A _ \mu (x))^n]\le b^n E[(\tau _ \mu  (x))^n].$

\subsection{The exit probability of $\mathcal X (t)$ from the left end of the interval $(0,b)$}
From the general theory of jump-diffusion processes with infinitesimal generator $\mathcal L$ (see e.g.  \cite{abundo:pms00}, \cite{gihman:1972} ) it is known that
the exit probability of ${\cal X} (t)$ from the left end of the interval $(0,b),$ when starting from $x \in (0,b),$ that is,
$ \pi _0 ^ \mu (x)= P [\mathcal X  (\tau _ \mu (x))=0 ],$  is the solution of the differential problem with boundary conditions:
\begin{equation}
\begin{cases}
{\cal L}u(x)=  0, \ x \in (0,b) \\
u(0) = 1, u(b) =0.
\end{cases}
\end{equation}
Thus, since we consider drifted BM with resetting, the generator ${\cal L}$ is given by \eqref{generator1}, and so $ \pi _0 ^ \mu (x)$ satisfies:
\begin{equation} \label{eqpi0}
\begin{cases}
{\cal L}u(x)= \frac 1 2  u''(x) + \mu  u'(x) - r  u (x) + r u (x_R) = 0, \ x \in (0,b) \\
u(0) = 1, u(b) =0.
\end{cases}
\end{equation}
Of course, the exit probability from the right end:
\begin{equation}
\pi _b (x)= P [\mathcal X  (\tau _ \mu (x))=b] = 1 -  \pi _0 (x)
\end{equation}
satisfies the same problem with boundary conditions $u(0) = 0, u(b) =1.$
The explicit solution of Eq. \eqref{eqpi0} is (see e.g. Eq. (25) of \cite{guoyan:24}):
\begin{equation} \label{pi0driftBMreset}
\pi _0 ^ \mu (x) = \frac {\sinh \left ((b-x_R)\sqrt {\mu ^2 +2r} \right ) + e^ {\mu (b-x)} \sinh \left ((x_R -x)\sqrt {\mu ^2 +2r} \right )} {\sinh \left ((b-x_R)\sqrt {\mu ^2 +2r} \right ) + e^ {b \mu} \sinh \left (x_R\sqrt {\mu ^2 +2r} \right )},
\end{equation}
while the probability of exit from the right end of the interval $(0,b)$ is:
\begin{equation}
\pi _b ^ \mu (x) = 1 - \pi _0 (x) = \frac {e^ {b \mu} \sinh \left (x_R\sqrt {\mu ^2 +2r} \right ) - e^ { \mu (b-x)} \sinh \left ((x_R -x)\sqrt {\mu ^2 +2r} \right )}{\sinh \left ((b-x_R)\sqrt {\mu ^2 +2r} \right ) + e^ {b \mu} \sinh \left (x_R\sqrt {\mu ^2 +2r} \right )}.
\end{equation}
\indent For $\mu =0,$ one gets the exit probability of (undrifted) BM with resetting:
\begin{equation} \label{pi0BMreset}
\pi ^0_0 (x) =  \frac { \sinh ((b-x_R) \sqrt {2r}) + \sinh ((x_R -x) \sqrt {2r}))} {\sinh ((b-x_R) \sqrt {2r}) + \sinh (x_R  \sqrt {2r}))}.
\end{equation}
If $x = x_R, $ then:
\begin{equation} \label{pi0xRBMreset}
\pi _0 ^ \mu (x_R) = \frac {\sinh \left ((b-x_R)\sqrt {\mu ^2 +2r} \right ) } {\sinh \left ((b-x_R)\sqrt {\mu ^2 +2r} \right ) + e^ {b \mu} \sinh \left (x_R\sqrt {\mu ^2 +2r} \right )}.
\end{equation}
Letting $r$ go to infinity in \eqref{pi0xRBMreset}, one obtains
\begin{equation}
\lim _ {r \rightarrow + \infty} \pi _0 ^ \mu (x) =
\begin{cases} \label{exitprorinfinity}
1, & x_R \in (0, b/2) \\
\frac 1 {1 + e^ {\mu b}}, & x_R= b/2 \\
0, & x_R \in (b/2, b).
\end{cases}
\end{equation}

Taking the limit as $r \rightarrow 0 ^+$ in \eqref{pi0driftBMreset},
one gets the well-known result for drifted  BM (without resetting), that is (see e.g. \cite{karlin2}, pg. 205 or \cite{borodin:1996}, pg. 233):
\begin{equation} \label{pidriftBMr=0}
\lim _ {r \rightarrow 0^+} \pi ^\mu _0(x)= \frac {e^ {2(b-x)\mu}- 1} {e^ {2b \mu}-1}, \ x \in (0,b).
\end{equation}
Finally, letting  $\mu$ go to zero in \eqref{pidriftBMr=0},  one obtains the exit probability of (undrifted) BM without resetting, i.e.
\begin{equation}
\lim _ {r \rightarrow 0^+} \pi _0 ^0(x) = \frac 1 b (b-x).
\end{equation}

\begin{figure}
\centering
\includegraphics[height=0.39 \textheight]{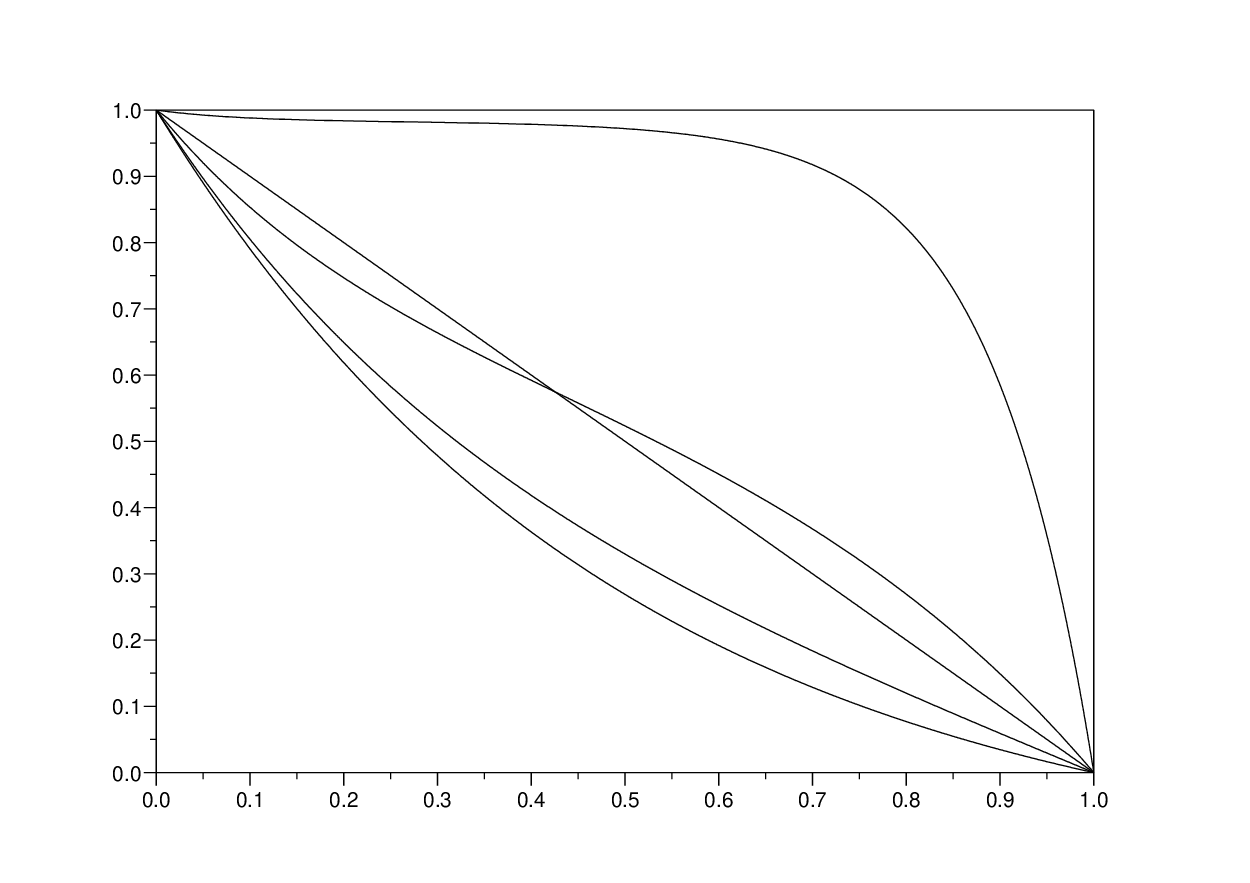}
\caption{Five examples of the graphs of $\pi_0^\mu (x)$ given by \eqref{pi0driftBMreset}, as function of $x \in (0,b),$ for
$b=1, \ \mu =1, \ r = 1, \ x_R= 1/4$ (lowest curve); $b=1, \ \mu =1,  r = 0$ (slightly higher curve); $b=1, \ \mu =0, \ r = 0 $ (straight line); $b=1, \ \mu =1, \ r = 5, \ x_R= 1/4$ (curve with inflection point);
\ $b=1, \ \mu =1, \ r = 50, \ x_R= 1/4$ (highest, concave curve).
}
\label{exitpro}
\end{figure}

\begin{Remark} \label{remark1}
From \eqref{pi0driftBMreset} we can  easily obtain the exit probability of $\mathcal X(t)$ from the left end of any interval $(a,b), \ a<b, $ under the condition that $\mathcal X(0)=x \in (a,b).$ In fact, let
$\pi ^- _{(a,b)}(x; \mu, x_R, r)$
denote the exit probability from the left end of the interval $(a,b),$ under the condition that
$\mathcal X(0) = x \in (a,b),$ in particular, with the previous notation it holds $\pi ^- _{(0,b)}(x; \mu, x_R, r) = \pi _0 ^ \mu (x);$ then, it is easy to see that
\begin{equation}
\pi ^- _{(a,b)}(x; \mu, x_R, r) = \pi ^- _{(0,b-a)}(x-a; \mu, x_R-a, r).
\end{equation}
Thus, $\pi ^- _{(a,b)}(x; \mu, x_R, r)$ can be obtained from \eqref{pi0driftBMreset}, by replacing $b$ with  $(b-a), \ x$  with $x-a,$ and $x_R$ with $x_R-a.$
\end{Remark}
In the Figure \ref{exitpro} we report five examples of the graph of $\pi_0^\mu (x)$ given by \eqref{pi0driftBMreset},  as a function of $x \in (0,b),$  for
$b=1, \ \mu =1, \ r = 1, \ x_R= 1/4$ (lowest curve); $b=1, \ \mu =1,  r = 0$ (slightly higher curve); $b=1, \ \mu =0, \ r = 0 $ (straight line); $b=1, \ \mu =1, \ r = 5, \ x_R= 1/4$ (curve with inflection point);
\ $b=1, \ \mu =1, \ r = 50, \ x_R= 1/4$ (highest, concave curve). Note that for $r=50$ (large resetting rate),
and starting point $x \in (0, 1/4),$
the exit probability  from the left end of $(0,b)$ is approximately one, in agreement with the theoretical result \eqref{exitprorinfinity}.

\newpage

\subsection{The Laplace transforms of the survival function and of the first-exit time}
Let $Q^ \mu _r(x,t;x_R)= P[ \tau _ \mu(x)>t]$ be the survival function of the FET $\tau _ \mu (x),$ where we use the notation of \cite{guoyan:24}, \cite{huang:24}, which makes explicit the dependence on $r$ and $x_R.$ Then, $Q^ \mu _r(x,t;x_R)$ satisfies the backward equation (see e.g. \cite{guoyan:24}, \cite{huang:24}):
$$
\frac {\partial} {\partial t} Q^ \mu _r(x,t; x_R) = {\cal L} Q ^ \mu _r(x,t;x_R)
$$
\begin{equation} \label{eqsurvivalprob}
= \frac 1 2 \frac {\partial ^2} {\partial x^2 }Q ^ \mu _r(x,t; x_R) + \mu \frac {\partial} {\partial x }Q ^ \mu _r(x,t; x_R) +r \left [Q ^ \mu _r(x_R,t; x_R) -Q^ \mu _r(x,t; x_R) \right ]
\end{equation}
with the boundary conditions
\begin{equation}
 Q ^ \mu _r(0,t; x_R) = Q^ \mu _r(b,t; x_R)=0,
 \end{equation}
 and the initial condition
 \begin{equation}
 Q^ \mu _r(x,0;x_R) = 1 .
 \end{equation}
The operator ${\cal L},$ defined by  \eqref{generator1}, is the infinitesimal generator of drifted BM with resetting $\mathcal X(t);$ it is meant as an operator which acts on  $Q ^ \mu _r(0,t; x_R),$ as a function of $x.$
For a general diffusion with resetting, one should take instead the operator given by \eqref{generator}. \par\noindent
Note that, by using a renewal argument (see \cite{guoyan:24}), one obtains:
\begin{equation} \label{renewalforQ}
Q^ \mu _r(x,t;x_R) = e^ {-rt}Q ^ \mu _0(x,t)+r \int _0 ^t ds e^ {-rs }Q ^ \mu _0(x_R,s)Q^ \mu _r(x,t-s;x_R),
\end{equation}
where $Q^ \mu_0(x,t)$ is the survival function of $\tau _ \mu (x)$ for the process without resetting, i.e. $r=0.$ \par\noindent
For $\lambda >0,$ let $\widehat Q^ \mu _r(x, \lambda; x_R)= \int_0^ \infty e^ {-\lambda t } Q ^ \mu _r(x,t; x_R) dt$ be the Laplace transform (LT) of
$Q^ \mu _r(x,t;x_R),$ and  let
$\widehat Q^ \mu _0(x, \lambda )= \int_0^ \infty e^ {-\lambda t } Q^ \mu _0(x,t) dt$ the LT of
$Q^ \mu _0(x,t).$
Then from \eqref{renewalforQ} one obtains:
\begin{equation} \label{renewalforQ2}
\widehat Q^ \mu _r(x, \lambda; x_R) = \frac {\widehat Q^ \mu _0(x,\lambda +r)} {1-r\widehat Q^ \mu _0(x_R,\lambda +r) }.
\end{equation}
By taking the Laplace transform in both sides of \eqref{eqsurvivalprob}, we get that $\widehat Q^ \mu _r(x, \lambda; x_R)$
satisfies, for $\lambda >0,$  the equation:
 \begin{equation} \label{eqLTsurvivalprob}
\frac 1 2 \frac {\partial ^2} {\partial x^2 }\widehat Q^ \mu_r(x, \lambda; x_R) + \mu  \frac {\partial} {\partial x }Q^ \mu _r(x, \lambda; x_R) - (r + \lambda) \widehat Q^ \mu _r(x, \lambda; x_R) + r \widehat Q^ \mu _r(x_R, \lambda; x_R) = -1,
 \end{equation}
 subject to the boundary conditions
  \begin{equation}
\widehat Q^ \mu _r(0, \lambda; x_R) = \widehat Q^ \mu _r(b, \lambda; x_R)=0.
 \end{equation}
Its explicit solution is (see e.g. Eq. (14) of \cite{guoyan:24}, with $D=1/2 ):$
\begin{equation} \label{explicit LTQdriftBMreset}
\widehat Q^ \mu _r(x, \lambda; x_R) = \frac {e^{-b \mu}\sinh(b \beta _ \lambda) - e^ {- (b+x) \mu }\sinh(\beta _ \lambda (b-x)) - e^ {- \mu x}\sinh (\beta _ \lambda x) }
{\lambda e^{- b \mu} \sinh (\beta _ \lambda b ) + r e^ {-(b+x_R) \mu}\sinh (\beta _ \lambda (b-x_R)) + r e^ {-x_R \mu }\sinh (\beta _ \lambda x_R )},
\end{equation}
where $\beta _ \lambda = \sqrt {\mu ^2 + 2(\lambda +r)}$ (we omit to explicit the dependence on $\mu $ and $r,$ for simplicity). \par\noindent
We denote by $M _ \mu (x, \lambda)$ the LT of $\tau _ \mu (x),$ that is,
 $$M _ \mu (x, \lambda) =  E \left [e^ {- \lambda  \tau(x)} \right ] = \int_0^ \infty e^ {-\lambda t } f_ {\tau _ \mu} (t) dt,$$
 where $f_ {\tau _ \mu} (t) = \frac d {dt} P[ \tau _ \mu(x) \le t ] $ is the probability density of $\tau _ \mu(x)$ (we omit to explicit the dependence of $M _ \mu (x, \lambda)$ on $r$ and $x_R).$
We have:
  \begin{equation} \label{QversusM}
\widehat Q^ \mu _r(x, \lambda;x_R) = \frac 1 \lambda \left [ 1 - M _ \mu  (x, \lambda) \right ],
 \end{equation}
or
\begin{equation} \label{MversusQ}
M _ \mu (x, \lambda) = 1 - \lambda \widehat Q ^ \mu _r(x, \lambda;x_R) .
 \end{equation}
In fact:
$$
\widehat Q^ \mu _r(x, \lambda; x_R)= \int_0^ \infty e^ {-\lambda t } Q^ \mu _r(x,t; x_R) dt = \int_0^ \infty e^ {-\lambda t } \left [ 1- P( \tau _ \mu(x) \le t \right ] dt $$
\begin{equation} \label{piece}
= \frac 1 \lambda - \int_0^ \infty e^ {-\lambda t } P(\tau _ \mu(x) \le t ) dt.
\end{equation}
Hence:
$$ M _ \mu (x, \lambda) =   \int _0 ^ \infty e^ {-\lambda t } f_ {\tau  _\mu }(t) dt = \int_0^ \infty e^ {-\lambda t } \frac d {dt} P [\tau _ \mu (x) \le t ] dt =$$
(integrating by parts)
\begin{equation}
= \lambda \int_0^ \infty e^ {-\lambda t } P [\tau _ \mu (x) \le t ] dt .
\end{equation}
The last integral is equal to $1/ \lambda - \widehat Q^ \mu _r(x, \lambda; x_R),$ tanking to \eqref{piece};
then \eqref{QversusM} and \eqref{MversusQ} soon follow.
By using \eqref{QversusM}, one obtains that $M _ \mu (x, \lambda),$ as a function of $x ,$ satisfies the differential problem (cf. Eq. (2.3) of \cite{abundo:FPA2023}):
 \begin{equation} \label{eqLTtau}
{\cal L}M _ \mu (x, \lambda) = \frac 1 2  \frac {\partial ^2} {\partial x^2 }M _\mu (x, \lambda) + \mu \frac {\partial} {\partial x } M _ \mu (x, \lambda)- (r + \lambda) M _ \mu (x, \lambda) + r M _ \mu  (x_R, \lambda) = 0, \ x \in (0,b)
 \end{equation}
with boundary conditions
  \begin{equation}
M _ \mu (0, \lambda)= M _ \mu  (b, \lambda) =1.
 \end{equation}
From \eqref{explicit LTQdriftBMreset}
the explicit expression of  $M_ \mu  (x, \lambda) = 1 - \lambda \widehat Q^ \mu _r (x, \lambda ; x_R)$ soon follows:
\begin{equation} \label{LTtaudriftBMreset}
M_ \mu (x, \lambda) =  1 - \frac {\lambda [e^{-b \mu}\sinh(b \beta _ \lambda) - e^ {- (b+x) \mu }\sinh(\beta _ \lambda (b-x)) - e^ {- \mu x}\sinh (\beta _ \lambda x) ]}
{\lambda e^{- b \mu} \sinh (\beta _ \lambda b ) + r e^ {-(b+x_R) \mu}\sinh (\beta _ \lambda (b-x_R)) + r e^ {-x_R \mu } \sinh (\beta _ \lambda x_R )}.
\end{equation}
In particular:
\begin{equation} \label{LTtauxRdriftedBM}
M _ \mu (x_R, \lambda) =  \frac {(\lambda +r) e^{-\mu x_R} \left [\sinh (\beta_ \lambda x_R) + e^ {-b \mu } \sinh (\beta_ \lambda (b-x_R))\right ] }
{ \lambda e^ {- b \mu} \sinh( \beta_ \lambda b) +
r e^ {-(b+x_R) \mu}\sinh (\beta _ \lambda (b-x_R)) + r e^ {-x_R \mu } \sinh (\beta _ \lambda x_R )},
\end{equation}
that, for $r=0,$  coincides with formula 3.01 of \cite{borodin:1996}, pg. 233.
\bigskip

For $\mu = 0,$ one gets from \eqref{explicit LTQdriftBMreset}:
\begin{equation} \label{explicit LTQBMreset}
\widehat Q^ 0_r(x, \lambda; x_R) = \frac {\sinh(\alpha _ \lambda b) - \sinh(\alpha _ \lambda x) - \sinh (\alpha _ \lambda (b-x)) } {\lambda  \sinh (\alpha _ \lambda b ) + r \sinh (\alpha _ \lambda x_R ) + r \sinh (\alpha _ \lambda (b-x_R) )},
\end{equation}
and from \eqref{LTtaudriftBMreset}:
\begin{equation} \label{LTtauBMreset}
M _ 0 (x, \lambda) =  \frac { r \sinh ( \alpha _ \lambda x_R ) + r \sinh (\alpha _ \lambda (b-x_R) ) + \lambda  \sinh ( \alpha _ \lambda x ) + \lambda  \sinh (\alpha _ \lambda (b-x) ) ] }{ \lambda  \sinh (\alpha _ \lambda b ) + r \sinh (\alpha _ \lambda x_R ) + r \sinh (\alpha _ \lambda (b-x_R) ) },
\end{equation}
where
$\alpha _ \lambda = \sqrt {2(\lambda +r)} $ (we omit the dependence on  $r).$
\par\noindent
In particular, for $x= x_R:$
\begin{equation} \label{MforxR}
M _ 0 (x_R, \lambda) = \frac {  \left (r + \lambda \right ) \left [ \sinh ( \alpha _ \lambda x_R ) +  \sinh (\alpha _ \lambda (b-x_R) ) \right ] }{ \lambda  \sinh (\alpha _ \lambda b ) + r \sinh (\alpha _ \lambda x_R ) + r \sinh (\alpha _ \lambda (b-x_R) ) } .
\end{equation}
Letting $b$ go to $+ \infty $ in \eqref{MforxR}, one obtains
$$ \lim _ {b \rightarrow + \infty} M _0  (x_R, \lambda) = \frac {(\lambda +r ) e^ {- x_R \sqrt {2(r + \lambda)}}} {\lambda + r e^ {- x_R \sqrt {2(r + \lambda)}} } $$
which coincides with Eq. (3.4) of \cite{abundo:FPA2023}.
\bigskip

For $\mu = r=0,$ namely, for (undrifted) BM without resetting, we obtain
\begin{equation}
\widehat Q^0_0(x, \lambda; x_R) = \frac {\sinh(\alpha _ \lambda b) - \sinh(\alpha _ \lambda x) - \sinh (\alpha _ \lambda (b-x)) } {\lambda  \sinh (\alpha _ \lambda b )  },
\end{equation}
and
\begin{equation}
M _ 0 (x, \lambda) = E [e^ {- \lambda \tau _ 0 (x)}] = 1 - \lambda \widehat Q^ 0 _0 (x, \lambda; x_R) = \frac { \sinh(\alpha _ \lambda x) + \sinh (\alpha _ \lambda (b-x)) } { \sinh (\alpha _ \lambda b )  },
\end{equation}
that, after straightforward algebraic manipulations, coincides with the well-known Darling's and Siegert's result \cite{darling:1953} (see also \cite{abundo:stapro13}):
\begin{equation} \label{LTDS}
M _ 0 (x, \lambda) =E [e^ {- \lambda \tau _ 0(x)}] = \frac {\cosh \left (\sqrt {2 \lambda} \ (x - \frac b 2 ) \right )} {\cosh \left (\sqrt {2 \lambda} \  \frac b 2 \right )}.
\end{equation}
Note that in all cases the LT of the FET, as a function of $\lambda ,$ turns out to be  well defined and finite even  for $\lambda$ belonging to a left neighborhood of $\lambda = 0.$ This is also true for $\mu = r=0$ (see
\eqref{LTDS}), as follows by using Euler's formula. Thus, all the moments of the FET are finite, and therefore also the moments of the FEA are finite.

\subsection{The moments of the FET} \par\noindent
Since the LT of $\tau _ \mu (x)$
is finite for $\lambda$ belonging to a neighborhood of $0,$
then the $n-$th order moments of $\tau _ \mu (x)$
exist  finite, and they are given by:
\begin{equation} \label{ODE1forM}
T_ n (x) := E \left [  \left ( \tau _ \mu (x)
\right ) ^n \right ]  = (-1)^n \left [ \frac {\partial
^n } {\partial \lambda ^n } M ^ \mu  (x, \lambda ) \right ] _ { \lambda =0} \ , n=1, 2 , \dots .
\end{equation}
The moments of $\tau _ \mu (x)$ can be also obtained from the LT $\widehat Q^ \mu _r(x, \lambda; x_R)$ of the survival function of $\tau _ \mu(x).$
As easily seen,  the first and second moments of $\tau _ \mu(x)$ are  given by:
\begin{equation} \label{momentsoftauasQ}
E[\tau(x)] = \widehat Q^ \mu _r(x, 0; x_R), \ \ E[\tau _ \mu ^2(x)]= - 2 \frac \partial {\partial \lambda} \widehat Q^ \mu _r(x, \lambda; x_R) |_ {\lambda =0}.
\end{equation}
By setting $T_0(x)=1$ and calculating  the $n-$th derivative with respect to $\lambda ,$ at $\lambda =0,$ of
both members  of   \eqref{eqLTtau}, we also obtain that
the $n-$th order moments $T_n(x)$  satisfy the ODEs (see e.g. \cite{abundo:TPMS2024} or \cite{abundo:pms00}):
\begin{equation} \label{eqmoments}
Lu(x) = \frac 1 2  u''(x) + \mu u'(x) = -n T_{n-1} (x) + r u(x) - r u(x_R) , \ x \in (0,b),
\end{equation}
with  the boundary conditions  $u(0)=u(b)= 0.$
Note that for $r=0$, Eq. \eqref{eqmoments} becomes  the celebrated Darling and Siegert's equation  (\cite{darling:1953}) for the moments of the first-passage time of a diffusion without resetting.
\bigskip

\noindent In particular, $E[\tau _ \mu (x)]= \widehat Q^ \mu _r (x,0; x_R),$
as a function of $x,$ is the solution of the differential equation with boundary conditions:
\begin{equation} \label{eqmeanFETBMreset}
\begin{cases}
\frac 1 2  u''(x) + \mu  u'(x) - r  u (x) + r u (x_R) = -1 , \ x \in (0.b) \\
u(0) = u(b) =0,
\end{cases}
\end{equation}
whilst $E[\tau _ \mu ^2(x)]$ is the solution of the differential problem:
\begin{equation} \label{eqsecondFETBMreset}
\begin{cases}
\frac 1 2  u''(x) + \mu  u'(x) - r  u (x) + r u (x_R) = -2 E[\tau (x)], \ x \in (0.b) \\
u(0) = u(b) =0,
\end{cases}
\end{equation}
To obtain the explicit expressions of $E[\tau _ \mu (x)]$ and $E[\tau _ \mu ^2(x)],$ instead of solving  \eqref{eqmeanFETBMreset} and \eqref{eqsecondFETBMreset}, is more convenient to use  formulae \eqref{momentsoftauasQ},
namely $E[\tau _ \mu (x)] = \widehat Q^ \mu _r(x, 0; x_R), \ E[\tau _ \mu ^2(x)]= - 2 \frac \partial {\partial \lambda} \widehat Q^ \mu _r(x, \lambda; x_R) |_ {\lambda =0},$
where $\widehat Q^ \mu _r(x, \lambda; x_R)$ is given by \eqref{explicit LTQdriftBMreset}.
Alternatively, $E[\tau _ \mu (x)]$ and $E[\tau _ \mu  ^2(x)]$ are respectively given by minus the derivative with respect to $\lambda$ of $M _ \mu (x, \lambda),$
and its second derivative, both calculated at $\lambda =0 .$ \par
\begin{figure}
\centering
\includegraphics[height=0.39 \textheight]{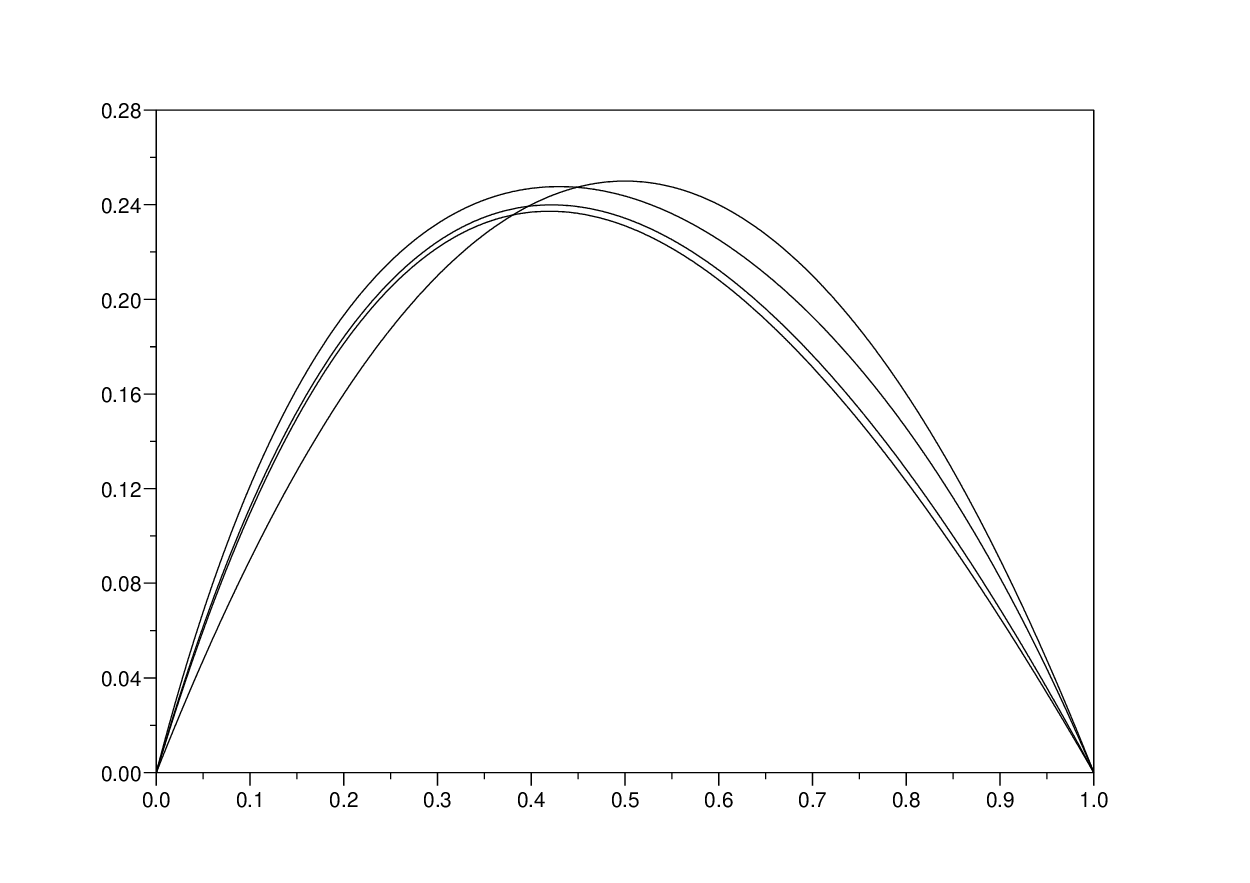}
\caption{Four examples of  plots of $E[\tau _ \mu (x)],$ as a function of $x \in (0,b);$
 the curves are ordered from the lower to the higher peak height. The sets of parameters are:
$b=1, \ \mu =1, \ r = 0$ (curve 1); $b=1, \ \mu =1, \ r = 1, \ x_R=1/4$ (curve 2);  $b=1, \ \mu =1, \ r = 5, \ x_R=1/4$ (curve 3);
$b=1, \ \mu = 0, \ r = 0 $ (curve 4).
}
\label{expectedFET}
\end{figure}

\noindent
In any case, one gets:
$$ E[\tau _ \mu (x)] = $$
\begin{equation} \label{EtauBMdriftreset}
= \frac {e^ {-b \mu }\sinh \left (b \sqrt {\mu ^2 +2r} \ \right )- e^ {- \mu x}\sinh \left (x \sqrt {\mu ^2 +2r} \ \right ) - e^ {-(b+x) \mu}\sinh \left ((b -x) \sqrt {\mu ^2 +2r} \ \right )   }  {r \left [ e^ {-x_R \mu}\sinh \left (x_R \sqrt {\mu ^2 +2r} \ \right ) + e^ {-(b+x_R) \mu}\sinh \left ((b - x_R) \sqrt {\mu ^2 +2r}  \ \right ) \right ]},
\end{equation}
The graph of $E[\tau _ \mu (x)],$ as a function of $x \in (0,b),$ shows a pseudo-parabolic behavior, with the only point of maximum inside the interval $(0,b),$ which is similar to the graph of the expected FET of BM without resetting.
In the Figure \ref{expectedFET} we report four examples of  plots of  $E[\tau _ \mu (x)],$ given by \eqref{EtauBMdriftreset}; the curves are ordered from the lower to the higher peak height. The sets of parameters are:
$b=1, \ \mu =1, \ r = 0$ (curve 1); $b=1, \ \mu =1, \ r = 1, \ x_R=1/4$ (curve 2);  $b=1, \ \mu =1, \ r = 5, \ x_R=1/4$ (curve 3);
$b=1, \ \mu = 0, \ r = 0 $ (curve 4).
Note that for the first three curves the peak is attained approximately at $\bar x = 0.45 ,$ while the curve 4, that refers to (undrifted) BM without resetting, attains the maximum exactly at $x=1/2,$ being in this case $E[\tau _0(x)] = x(b-x).$

\bigskip
Moreover, it holds:
$$E[\tau _ \mu ^2(x)] =$$
  $$ = \frac {2 \left [e^ {- \mu b}\sinh(b\beta_0)- e^ {- \mu x}\sinh(x\beta_0)- e ^ {-(b+x) \mu}\sinh((b-x)\beta_0) \right ] } {r^2 \left [e^ {- x_R \mu} \sinh (x_R \beta _0 )+ e^ {-(b+x_R) \mu }\sinh ((b-x_R)\beta _0 ) \right ]^2}  $$
$$ \times \left [e^ {- \mu b}\sinh (b\beta _0) + \frac r {\beta _0}  \left (x_R e^ {- \mu x_R}\cosh (x_R \beta _0) +
(b -x_R) e^ {-(b+x_R) \mu}\cosh (x_R \beta _0) \right )  \right ] +
$$
\begin{equation} \label{Etau2BMdriftreset}
- \frac 2 {\beta _0} \ \frac {b e^ {- \mu b}\cosh ( b \beta _0) - x e^ {- \mu x}\cosh ( x \beta _0) - (b-x) e ^ {- (b+x)\mu }\cosh ( (b-x) \beta _0)} {r \left [e^ {- x_R \mu} \sinh (x_R \beta _0 )+ e^ {-(b+x_R) \mu }\sinh ((b-x_R)\beta _0 ) \right ]},
\end{equation}
where $ \beta _0 = \sqrt {\mu  ^2 + 2r}.$ \par
Taking the limit as $r \rightarrow 0 ^+$ in \eqref{EtauBMdriftreset} one obtains the  formula for drifted BM without resetting:
\begin{equation} \label{EtaudriftBMnoreset}
\lim _ {r \rightarrow 0 ^+ } E [\tau _ \mu  (x)]= \frac 1 \mu \left [\frac {b(1- e^ {- 2 \mu x})} {1- e^ {-2\mu b}} -x \right ]
\end{equation}
(in the calculations, it is convenient to use the Taylor's expansion of $\sinh (\gamma \sqrt {\mu ^2 + 2r}),$ as $r \rightarrow 0 ).$

\begin{figure}
\centering
\includegraphics[height=0.39 \textheight]{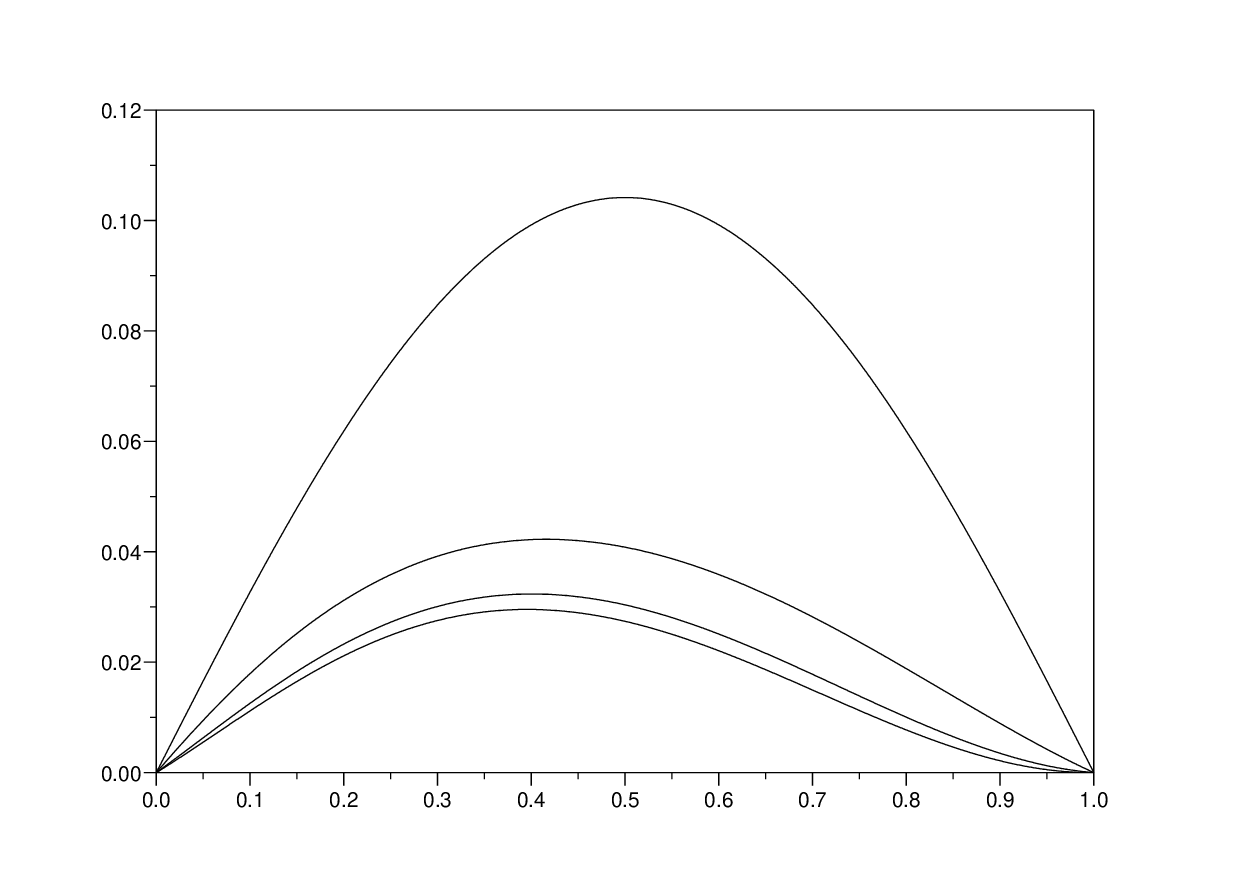}
\caption{Four examples of  plots of $E[\tau _ \mu ^2(x)],$ as a function of $x \in (0,b);$
 the curves are ordered from the lower to the higher peak height. The sets of parameters are the same ones as in the Figure
\ref{expectedFET}.
}
\label{expectsquareFET}
\end{figure}

In the Figure \ref{expectsquareFET} we report four examples of  plots of  $E[\tau _ \mu ^2(x)]$; the curves are ordered from the lower to the higher peak height.
The sets of parameters are the same ones as in the Figure \ref{expectedFET}.
Note that for the first three curve the peak is attained approximately at $\bar x = 0.45  ,$ while the curve 4, that refers to BM without resetting, attains the maximum exactly at $x=1/2.$
\bigskip

\bigskip

\bigskip

For $\mu =0,$ the corresponding expressions for the first two moments of $\tau_0 (x)$ are the same ones as those with drift, without the exponentials. Note that
\begin{equation}
\max _ {x \in (0,b)} E[\tau_0(x)]= E[\tau_0(b/2)],
\end{equation}
and, if $x_R \notin \{0, b \}:$
\begin{equation}
\lim _ {r \rightarrow + \infty }\max _ {x \in (0,b)} E[\tau_0(x)]= + \infty .
\end{equation}
\noindent For $x=x_R,$ one obtains:
\begin{equation} \label{EtauxR}
 E[ \tau _0 (x_R)]= \frac 1 r \left ( \frac {\sinh (b \sqrt{2r})} { \sinh (x_R \sqrt{2r})+ \sinh ((b-x_R) \sqrt{2r}) } -1 \right ),
\end{equation}
and letting $b$ go to $+ \infty ,$ one gets
\begin{equation}
 \lim _ {b \rightarrow + \infty} E[ \tau _0 (x_R)]= \frac 1 r \left (e^ {x_R \sqrt {2r}} -1  \right ) ,
\end{equation}
which coincides with Eq. (3.10) of \cite{abundo:FPA2023}, where is treated BM  with resetting with $b= + \infty.$
\bigskip

\noindent
As easily seen (by taking the limit in \eqref{EtaudriftBMnoreset} as $\mu \rightarrow 0),$ one finds
\begin{equation}
\lim _ {r \rightarrow 0^ +} E[\tau _0 (x)] = x(b-x), \ x \in (0,b),
\end{equation}
that is, the  well-known formula for the expected FET of (undrifted) BM without resetting. Moreover, it holds:
\begin{equation}
\lim _ {r \rightarrow 0^ +} E[\tau^2 _0 (x)] = \frac x 3 [b^2(2-b)+bx^2(x-2)], \ x \in (0,b),
\end{equation}
which is the expression of the second moment of the FET of (undrifted) BM without resetting.

\subsection{The Laplace transform of the first-exit area}
Recall that, for $x \in (0,b),$  the  first-exit area (FEA) is:
\begin{equation}
A _ \mu (x) = \int _0 ^ { \tau _ \mu (x)} {\cal X} (t) dt , \ x \in (0,b) .
\end{equation}
For $\lambda >0, $ let us denote by $N _ \mu  (x, \lambda)$  the LT  of $A _ \mu (x),$ that is,
\begin{equation}
N _ \mu  (x, \lambda)= E \left [e^ {- \lambda A_ \mu (x)} \right ].
\end{equation}
Then, $N_ \mu  (x, \lambda),$ as a function of $x \in (0,b),$  solves the following  differential problem with boundary conditions (see e.g. \cite{abundo:FPA2023}, or  \cite{abundo:MCAP2013}):
\begin{equation}
\begin{cases}
{\cal L} u(x) = \lambda x u(x) ,  \ x \in (0,b)\\
u(0) = u(b) =1 ,
\end{cases}
\end{equation}
that is,
\begin{equation} \label{ODEforNBMreset}
\begin{cases}
\frac 1 2 u''(x)+ \mu u'(x)  - \lambda x u(x) - r  u(x) + r u(x_R) =0, \ x \in (0,b) \\
u(0) = u(b) =1 .
\end{cases}
\end{equation}
Unfortunately, this is a differential equation with non-constant coefficients and it cannot be solved in terms of elementary functions, but only  special functions.\par\noindent
In the special case $\mu = r =0,$ Eq. \eqref{ODEforNBMreset} becomes $u''(x) - 2 \lambda x u(x) =0;$
two fundamental solutions are $u_1(x)= Ai \left ((2 \lambda ) ^ {1/3}x \right ),$ and $u_2(x)=Bi \left ((2 \lambda ) ^ {1/3}x \right ),$
where $Ai$ and $Bi$ are the
first and second kind Airy functions (see e.g. \cite{Abramowitz}).
Then, the solution can be written as $u(x)=  \alpha_1 u_1(x)+ \alpha_2 u_2(x),$
where
$ \alpha_1, \alpha_2 $ are constants with respect to $x$, to be determined (really, they depend on $b$ and $\lambda )$; by imposing the boundary conditions $u(0) = 1, u(b) =1,$ one finally finds that the solution of
\eqref{ODEforNBMreset} is (for $\mu = r =0):$
\begin{equation}
N _ 0 (x, \lambda) = \alpha_1 Ai \left ((2 \lambda ) ^ {1/3}x \right )+ \alpha_2 Bi \left ((2 \lambda ) ^ {1/3}x \right ),
\end{equation}
with:
\begin{equation}
\alpha _1 =  \frac {Bi(0) - Bi( (2 \lambda ) ^ {1/3}b)} {Bi (0) Ai ( (2 \lambda ) ^ {1/3}b)- Bi((2 \lambda ) ^ {1/3}b ) Ai(0) }, \
\alpha _2 = 1- \frac {\alpha _1 Ai(0)}{Bi(0)}.
\end{equation}
In the case $b= + \infty ,$ the analogous
formula for $N _ 0 (x, \lambda)$ with $r=0$ was found in \cite{singh} (see also \cite{abundo:FPA2023}).

\subsection{The moments of the FEA} \par\noindent
Since  $\tau _\mu (x)$ possesses finite moments of any order $n,$ also
the $n-$th order moments of $A_ \mu (x)$
exist  finite, and they are given by:
\begin{equation}
S_ n (x) := E \left [  \left ( A_ \mu (x)
\right ) ^n \right ]  = (-1)^n \left [ \frac {\partial
^n } {\partial \lambda ^n } N_ \mu  (x, \lambda) \right ] _ { \lambda =0} \ , n=1, 2 , \dots .
\end{equation}
By setting $S_0(x)=1$ and calculating  the $n-$th derivative with respect to $\lambda ,$ at $\lambda =0,$ of
both members  of Eq. \eqref{ODEforNBMreset}, we obtains that
the $n-$th order moments of $A_ \mu (x),$ i.e. $S_n(x),$  satisfy the ODEs (see e.g. \cite{abundo:TPMS2024} or \cite{abundo:MCAP2013}):
\begin{equation} \label{eqmomentsarea}
LS_n(x) = \frac 1 2 S_n''(x) + \mu S_n'(x)) -n x S_{n-1} (x) + r S_n(x) - r S_n(x_R) , \ x \in (0,b),
\end{equation}
with  the boundary conditions  $S_n(0)=S_n(b)= 0.$ \par
Now, we will calculate the first two moments of $A_ \mu (x).$
Because of complexity of calculations, we will limit ourselves to the case when $\mu =0,$ namely (undrifted) BM with resetting.
Note that, for $\mu =0, \ r=0$ (undrifted BM without resetting) and $b= + \infty,$ the moments of $A_0(x)$ are infinite (see \cite{abundodelvescovo:mcap17}); instead, for
$ \mu =0, \ r \ge 0$ and $b < + \infty ,$ the moments of  $A_ 0(x)$ are finite, because the moments of $\tau _ 0 (x)$ are finite.
\par
The Eq. \eqref{eqmomentsarea} for $n=1,$ which provides
 $E[A_0(x)],$ becomes:
\begin{equation} \label{eqmeanFEABMreset}
\begin{cases}
\frac 1 2 u''(x)  - r u(x) + r u(x_R) = -x , \ x \in (0,b)\\
u(0) = u(b) =0.
\end{cases}
\end{equation}
By solving it with standard methods, one finds:
\begin{equation} \label{meanFEA}
E[A _0(x)] =  c_1 e^{- x \sqrt {2r}} + c_2 e^{ x\sqrt {2r}} + \frac x r + c_3,
\end{equation}
where
\begin{equation} \label{costantscalpha}
\begin{cases}
c_1= \frac {b/r - c_3 (e^ {b \sqrt {2r}} -1)} {2 \sinh (b \sqrt {2r})}
 \\
c_2= \frac {-b/r - c_3 (1 - e^ {-b \sqrt {2r}})} {2 \sinh (b \sqrt {2r})}
\\
c_3 = \frac 1 r \ \frac {x_R \sinh (b \sqrt {2r}) -b \sinh (x_R \sqrt {2r}) } {\sinh ((b-x_R) \sqrt {2r}) + \sinh (x_R \sqrt {2r})}
\end{cases}
\end{equation}
(note that the constants $c_i$ depend on $b, \ r$ and $x_R ).$  \par
\bigskip

\noindent For $r=0$ Eq. \eqref{eqmeanFEABMreset} simply reads:
\begin{equation} \label{eqmeanFEABMr=0reset}
\begin{cases}
\frac 1 2 u''(x)  = -x , \ x\in (0,b)\\
u(0) = u(b) =0,
\end{cases}
\end{equation}
whose solution provides:
\begin{equation} \label{meanFEAr=0}
E[A _0(x)] =  \frac x 3 (b^2 - x^2) \ (r=0).
\end{equation}
\bigskip

The second order moment, $ E[A_0^2(x)],$ satisfies the differential problem:
\begin{equation} \label{eqsquareFEABMreset}
\begin{cases}
\frac 1 2 u''(x)  - r u(x) + r u(x_R) = -2x E[A_0(x)] , \ x \in (0,b)\\
u(0) = u(b) =0,
\end{cases}
\end{equation}
whose explicit solution is found by standard methods:
$$
E[A _0 ^2(x)] =   e^{- x \sqrt {2r}} \left [d_1 + \frac {c_1} {\sqrt {2r}} \left (x^2 + \frac {x}{\sqrt {2r}} \right ) \right ]+ e^{ x\sqrt {2r}}
\left [d_2 - \frac {c_2} {\sqrt {2r}} \left (x^2 - \frac {x}{\sqrt {2r}} \right ) \right ] +
$$
\begin{equation} \label{squareFEA}
+ \frac {2x^2} {r^2}+ \frac {2 c_3 x }r   + \frac {2} {r^3} + E[A _0 ^2(x_R)],
\end{equation}
where
$$
d_1= \Big \{ e^ {-b \sqrt {2r}} \frac {c_1} {\sqrt {2r}} \left (b^2 + \frac {b} {\sqrt {2r}} \right )- e^ {b \sqrt {2r}} \frac {c_2} {\sqrt {2r}} \left (b^2 - \frac {b} {\sqrt {2r}} \right )
- \frac 2 {r^3} \left (e^ {b \sqrt {2r}} -1 \right ) +
$$
\begin{equation}
- \left (e^ {b \sqrt {2r}} -1 \right ) E[A _0 ^2(x_R)] + \frac {2b^2} {r^2} + \frac {2 c_3 b} r \Big \} \frac 1 {2 \sinh (b \sqrt {2r}) } ,
\end{equation}
\begin{equation}
d_2 = - \frac {2} {r^3} - E[A _0 ^2(x_R)] - d_1,
\end{equation}
$$
E[A _0 ^2(x_R)]=
$$
$$
\Big [-2d_1 \sinh (x_R \sqrt {2r}) +
e^ {-x_R \sqrt {2r}} \frac {c_1} {\sqrt {2r}} \left (x_R^2 + \frac {x_R} {\sqrt {2r}} \right ) - e^ {x_R \sqrt {2r}} \frac {c_2} {\sqrt {2r}} \left (x_R^2 - \frac {x_R} {\sqrt {2r}} \right )
- \frac 2 {r^3}e^ {x_R \sqrt {2r}} +
$$
\begin{equation} \label{squareFEAxR}
+ \frac 2 {r^2} x_R ^2+ \frac {2 c_3 x_R} r + \frac 2 {r^3}
        \Big ] e^{-x_R \sqrt {2r}},
\end{equation}
and finally $c_i$  are given by \eqref{costantscalpha}.
\bigskip

\begin{figure}
\centering
\includegraphics[height=0.39 \textheight]{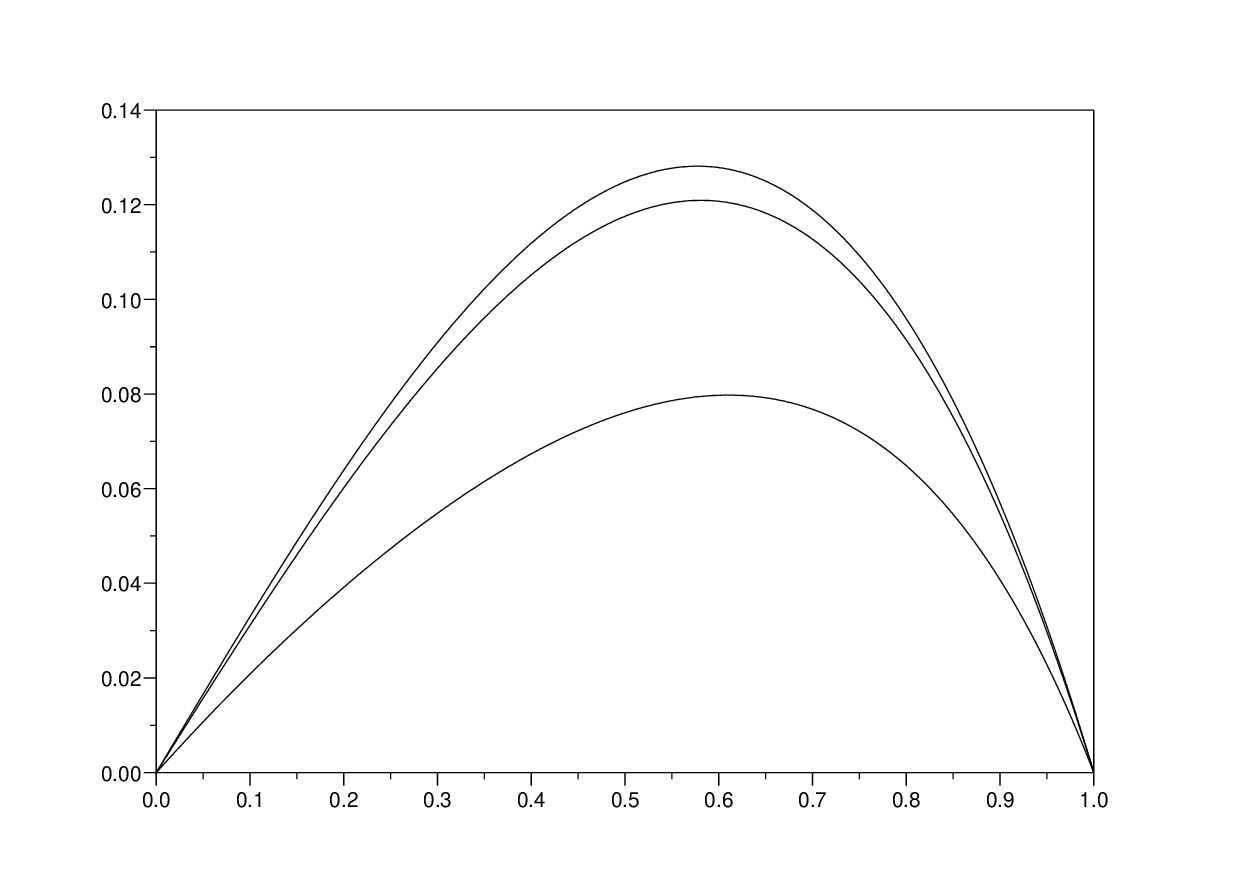}
\caption{Three examples of graphs of $E[A _ 0 (x)],$ as a function of $x \in (0,b);$
 the curves are ordered from the higher to the lower peak height. The sets of parameters are:
$\mu =0, \ b=1, r = 0$ (curve 1);
$\mu =0, \ b=1, r = 1/2, \ x_R= 1/8$ (curve 2); $\mu =0, \ b=1, r = 5, \ x_R= 1/8$ (curve 3).
}
\label{expectedFEAmu0}
\end{figure}


Note that to obtain the solutions to the differential problems \eqref{eqmeanFEABMreset} and \eqref{eqsquareFEABMreset} one has first to solve the associated homogeneous differential equation, and then find a
particular solution to be added to the solution of the homogeneous equation.
Actually, the calculations to obtain $E[A _0 (x)] $ and $E[A _0 ^2(x)]$ are far more complicated than those necessary to get the corresponding quantities, when $b= + \infty,$ since
in that case one of the constants $c_i$ and  $d_i$ vanishes (see \cite{abundo:FPA2023}). This is the reason why we decided to derive explicit formulae only in the simpler case when $\mu =0$  (undrifted BM with resetting). \par
\begin{figure}
\centering
\includegraphics[height=0.39 \textheight]{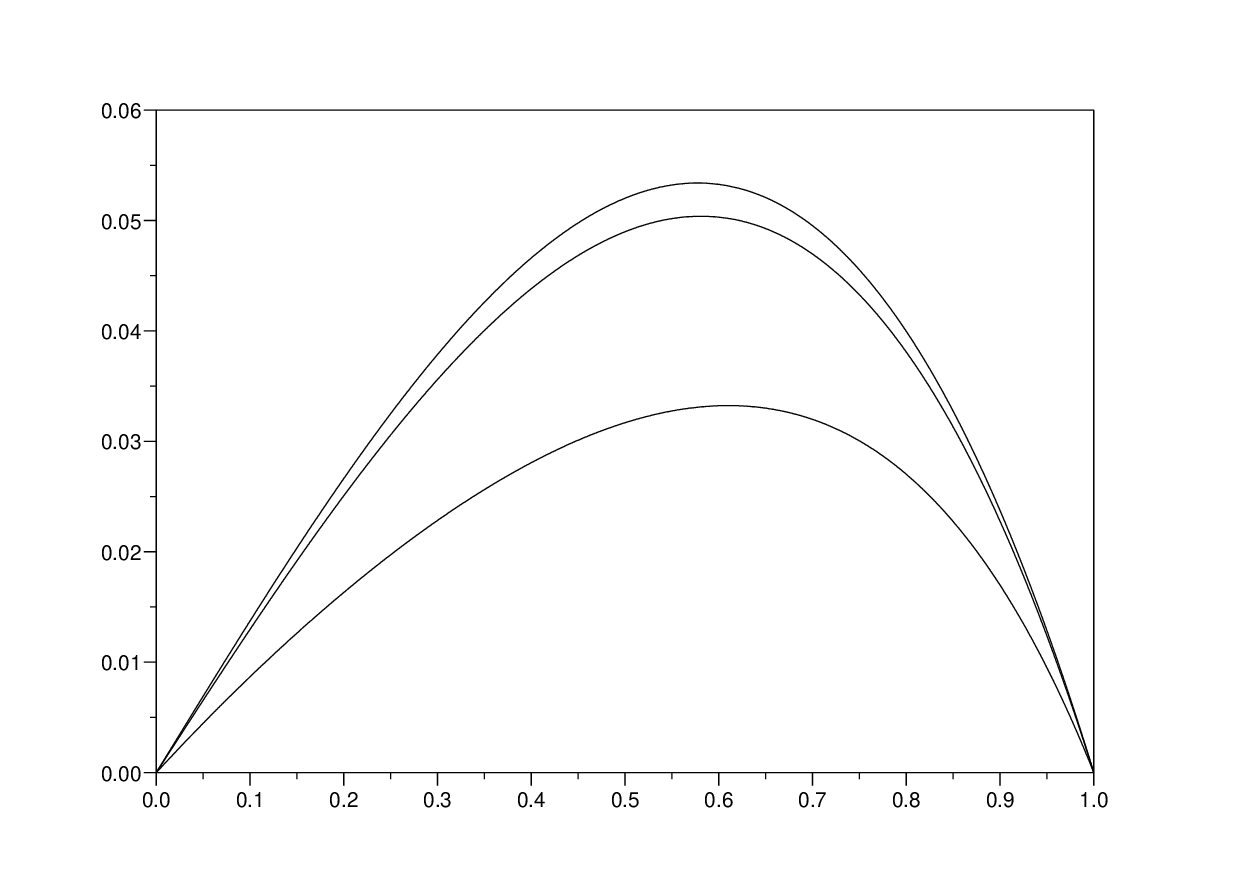}
\caption{Three examples of graphs of  of  $E[A ^2_ 0 (x)],$ as a function of $x \in (0,b);$
the curves are from the higher to the lower peak height. The sets of parameters are:
$\mu =0, \ b=1, r = 0$ (curve 1);
$\mu =0, \ b=1, r = 1/2, \ x_R= 1/8$ (curve 2); $\mu =0, \ b=1, r = 5, \ x_R= 1/8$ (curve 3).
}
\label{expectedsquareFEAmu0}
\end{figure}

\noindent
In the Figure \ref{expectedFEAmu0}, we report three examples of plots of $E[A _ 0 (x)],$ as a function of $x \in (0,b);$
the curves are ordered from the lower to the higher peak height. The sets of parameters are:
 $\mu =0, \ b=1, r =0$ (curve 1);
$\mu =0, \ b=1, r = 1/2, \ x_R= 1/8$ (curve 2); $\mu =0, \ b=1, r = 5, \ x_R= 1/8$ (curve 3).
In the Figure \ref{expectedsquareFEAmu0}, we report three examples of plots of $E[A _ 0 ^2(x)],$ as a function of $x \in (0,b);$
the curves are ordered from the lower to the higher peak height, and
the sets of parameters are  the same ones as in the Figure \ref{expectedFEAmu0}. \par\noindent
Note that for $r=0$ the expression  \eqref{squareFEA} for the second moment of $A_0(x)$ is not defined,
so one has to take the limit as $r$ tends to zero.

\subsection{Joint moment of $\tau_ \mu (x)$ and  $A_ \mu(x)$} \par\noindent
This subsection is devoted to find an  explicit expression of $E[\tau _ \mu (x) A_ \mu(x)],$ i.e. the joint moment of $\tau _ \mu(x)$ and  $A_ \mu(x).$
Actually, because of complexity of calculations, we will develop the calculations only in the case when $\mu =0.$ \par\noindent
The joint LT of $(\tau _ 0(x), A_ 0(x))$ is
$$
M_0 (x, \lambda _1, \lambda _2 ):= E \left [ e^ { - \lambda _1 \tau _0 (x)} e ^ { - \lambda _2 A_ 0(x)} \right ], \lambda _1, \lambda _2 >0 .
$$
By reasoning as in \cite{abundo:FPA2023} or \cite{abundo:MCAP2013}, we find that it satisfies the differential equation:
\begin{equation} \label{eqforjointLT}
(L- r) u(x) = ( \lambda _1 + \lambda _2 x) u (x) - r u (x_R), \ x \in (0,b),
\end{equation}
with the boundary conditions $u(0)= u(b)=0,$ where now $Lu(x)= \frac 1 2 u''(x),$ because $\mu =0.$ \par\noindent
Then, by taking
$\frac {\partial ^2} {\partial \lambda _1 \partial \lambda _2 } $ in both members of Eq. \eqref{eqforjointLT} and calculating it at $\lambda _1 = \lambda _2 =0,$ we obtain that
the joint moment $V(x):= E[\tau _0(x) A _ 0(x)]$  satisfies the boundary value differential problem:
\begin{equation} \label{eqforVxmu}
\begin{cases}
\frac 1 2 V''(x)  - r V(x) = -x E[ \tau _ 0(x)] -E[A _ 0(x)] -r V(x_R) , \ x \in (0,b)  \\
 V(0)= 0, \ V(b)=0.
\end{cases}
\end{equation}
Finding its explicit solution is very boring, however,
by standard methods, one gets:
\begin{equation} \label{explicitV}
 V(x) = (C_1(x) + \beta _1) e^ {-x \sqrt {2r}}+ (C_2(x) + \beta _2) e^ {x \sqrt {2r}} + V(x_R),
\end{equation}
where the functions $C_i(x), \ i=1,2 ,$ are:
$$
C_1(x)= \frac 1 {2 \sqrt {2r}} \Big \{ e^ { x \sqrt {2r}} \left [\left (\frac x {\sqrt {2r}} - \frac 1 {2r} \right ) \left ( \frac {D_2} {D_1}  + \frac 1 r \right ) + \frac {c_3} {\sqrt {2r}} \right ] $$
\begin{equation}
+  \frac {1 - e^ {b \sqrt {2r} }} {4 D_1}  x^2  + c_1 x+ e^ {2x \sqrt {2r}}  \left [ \frac {e ^ {-b \sqrt {2r}} -1} {8 D_1} \left ( \frac {2x} {\sqrt {2r}} - \frac 1 {2r} \right ) + \frac {c_2} {2\sqrt {2r}}\right ] \Big \}
\end{equation}
$$
C_2(x)= \frac 1 {2 \sqrt {2r}} \Big \{ e^ { -x \sqrt {2r}} \left [\left (\frac x {\sqrt {2r}} + \frac 1 {2r} \right ) \left ( \frac {D_2} {D_1}  + \frac 1 r \right ) + \frac {c_3} {\sqrt {2r}} \right ] +$$
\begin{equation}
  + \frac {1 - e^ {-b \sqrt {2r} }} {4 D_1}  x^2  - c_2 x+ e^ {-2x \sqrt {2r}}  \left [ \frac {1-e ^ {b \sqrt {2r}} } {2 D_1} \left ( \frac {x} {2\sqrt {2r}} + \frac 1 {8r} \right ) + \frac {c_1} {2 \sqrt {2r}}\right ] \Big \},
\end{equation}
with:
\begin{equation}
D_1 = r [\sinh(x_R \sqrt {2r}) + \sinh((b-x_R) \sqrt {2r}) ] , \ D_2 = \sinh(b \sqrt {2r}),
\end{equation}
and $c_1, \ c_2, \ c_3$ are given by \eqref{costantscalpha} (we omit the dependence of the functions $C_i(x)$ and the constants $D_i$ on $b, r, x_R).$
By imposing the boundary conditions $ V(0)= 0, \ V(b)=0,$ one finds the constants $\beta _i$ (actually, they depend on $b, r, x_R):$
\begin{equation}
\beta_1=  - \frac { V(x_R) (e^ {b \sqrt {2r}}-1) - C_1(b) e^ {-b \sqrt {2r}} - \left [ C_2(b) - C_1(0) - C_2(0) \right ] e^ {b \sqrt {2r}} } {2 \sinh (b \sqrt {2r})}
\end{equation}
\begin{equation}
\beta_2= - V(x_R) - C_1(0) - C_2(0) - \beta_1.
\end{equation}
\par\noindent
Finally, by taking $x=x_R$ in \eqref{explicitV}, one obtains an equation in the unknown $V(x_R),$ which allows to find it:
$$ V(x_R)= \Big \{ C_1(x_R) e^ {-2x_R \sqrt {2r}} + C_2(x_R) -C_1(b) \ \frac {\sinh(x_R \sqrt {2r})} {\sinh(b \sqrt {2r})} \ e^ {-(x_R +b) \sqrt {2r}}+$$
$$ -  \frac {\sinh(x_R \sqrt {2r})} {\sinh(b \sqrt {2r})} \ e^ {(b- x_R) \sqrt {2r}} \ \Big [ C_2 (b) - C_1(0) -C_2(0) \Big ] - C_1(0) -C_2(0) \Big \} $$
\begin{equation}
\times \Big \{ 1-  \frac {\sinh(x_R \sqrt {2r})} {\sinh(b \sqrt {2r})} \left (e^ {b \sqrt {2r}} -1 \right ) e^ {-x_R \sqrt {2r}} \ \Big \} ^ {-1}
\end{equation}
(note that the expression in $\{ \ \}$ is  positive for every value of $x_R ).$ \par\noindent
Thus, the explicit forms of
$E[\tau _0 (x) A _ 0(x)] = V(x)$ and of $Cov[\tau_0(x), A_0(x)] = E[\tau _0 (x) A _ 0(x)] - E[\tau _0 (x)] E[A _0 (x)]$ is obtained from \eqref{explicitV},
\eqref{EtauBMdriftreset} with $\mu =0,$  and \eqref{meanFEA}.
\par
\begin{figure}
\centering
\includegraphics[height=0.24 \textheight]{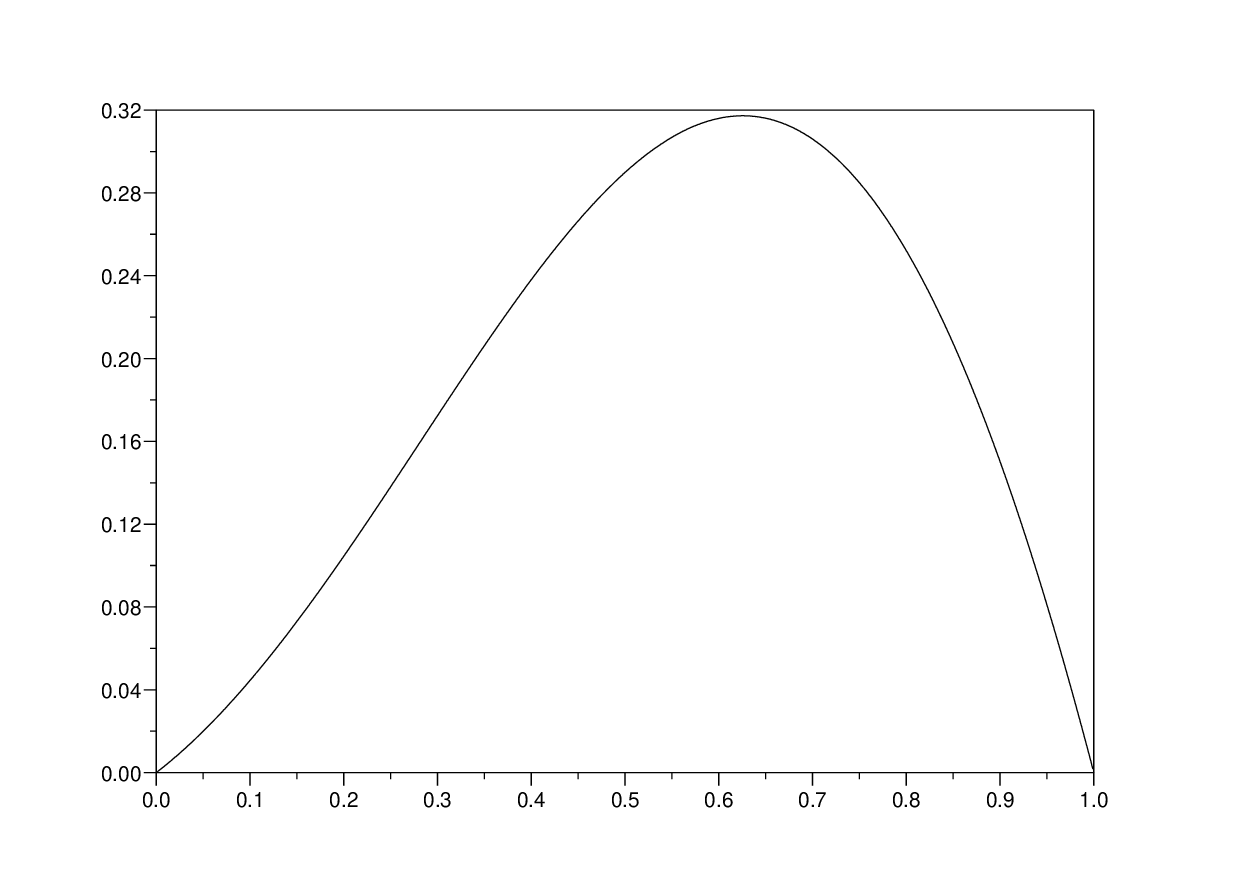}
\includegraphics[height=0.24 \textheight]{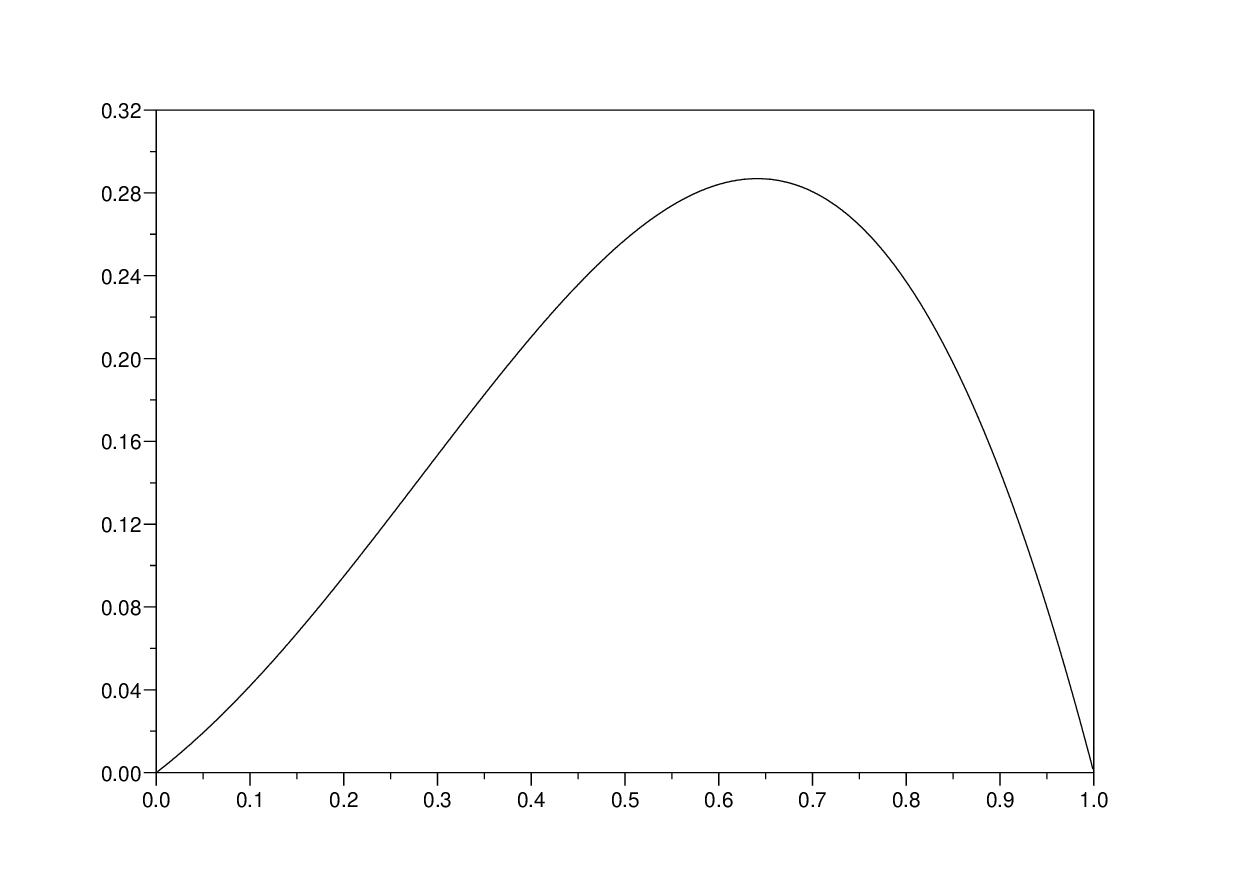}
\caption{Graphs of $E[\tau_0(x_R) A_0(x_R)]$ (left panel) and  $Cov [\tau_0(x_R), A_0(x_R) ]$ (right panel), as functions of $x_R \in (0,b),$ in the case of  BM with resetting, for
$b=1, r = 1/2.$
}
\label{EAtaucov}
\end{figure}

\noindent
In the Figure \ref{EAtaucov}, we report an example of the graph of $E[\tau _0 (x_R) A _ 0(x_R)]$ (left panel) and the graph of $Cov[\tau_0(x_R), A_0(x_R)]$ (right panel), as  functions of $x_R \in (0,b),$ for $b=1$ and  $r = 1/2.$ As we see, $Cov [\tau_0(x_R), A_0(x_R)]$ is positive, that is $\tau_0(x_R)$ and $A_0(x_R)$ are positively correlated; the covariance starts from zero, it reaches its maximum in the interior of $(0,b),$ and then it vanishes at $x_R =b$ again (compare this behavior with that of $Cov [\tau_0(x), A_0(x)]$ in the case of BM with resetting with $b = + \infty ,$  see  \cite{abundo:FPA2023}).
\bigskip
%
%

\subsection{The distributions of the maximum and minimum displacement till the exit time}
The maximum displacement of ${\cal X}(t),$ starting from $x \in (0,b),$  till the FET $\tau _ \mu (x)$, is the  r.v.  $\mathcal M _x = \max _ { t \in [0, \tau _ \mu  (x)]} {\cal X}(t)$
(of course one has $b \ge \mathcal M _x \ge x).$ Two cases are possible: in the first one ${\cal X}(t)$ first exits the interval $(0,b)$ through $0,$ and one can only say that
   $\mathcal M _x \ge x,$ in the second case
${\cal X}(t)$ first exits the interval $(0,b)$ through $b$ and $ \mathcal M _x =b.$
Let us denote by $E_0$ and $E_b,$ the event that ${\cal X}(t)$ first exits the interval $(0,b)$ through $0,$ respectively $b,$ namely
$E_0 = \{ \mathcal X (\tau _ \mu (x))= 0 \}, \ E_b = \{ \mathcal X (\tau _ \mu (x))= b \}.$
For $b\ge z \ge x,$
the event  $\{ \mathcal M _x \le z , \ \mathcal X (\tau _ \mu (x))= 0  \}$ occurs if and only if ${\cal X}(t)$ first
exits  the interval $(0, z)$ through the left end $0.$
Therefore, for $b \ge z \ge x$  one gets that  $w_ {\mathcal M}(z):=P(  \mathcal M _x \le z , \ {\cal X}(\tau _ \mu (x))= 0  )$ is the solution of
the differential boundary problem (see \cite{abundo:FPA2023}, \cite{abundo:MCAP2013}):
\begin{equation} \label{eqmaxdisplX}
\begin{cases}
 L u(x) + r u(x_R) - r u(x) =0, \ x \in (0,z) \\
u(0) =1, \ u(z) =0 ,
\end{cases}
\end{equation}
with  $L u = \frac 1 2 u''(x) + \mu u'(x).$
\par\noindent
Its explicit form is obtained from  \eqref{pi0driftBMreset}, by replacing the interval $(0,b)$ with $(0,z),$ so obtaining:
\begin{equation} \label{distrintM}
P \left [ \mathcal M_x \le z, {\cal X}(\tau _ \mu (x)) \right ] =0) =
 \begin{cases}
  0, & z < x \\
 \frac { \sinh \left ((z-x_R) \sqrt {\mu ^2 + 2r} \right ) + e^ {\mu (z-x)}\sinh \left ((x_R -x) \sqrt {\mu ^2 + 2r} \right  )} {\sinh \left ((z-x_R) \sqrt {\mu ^ 2 + 2r} \right ) + e^ {z \mu }\sinh \left (x_R  \sqrt {\mu ^2 + 2r} \right )}, &  x \le z \le b  .
\end{cases}
\end{equation}
Thus $ w_ {\mathcal M}(z)$ increases from $0$ to $\pi ^ \mu  _0 (x).$ \par\noindent
The conditional distribution function of $\mathcal M _x$ is:
\begin{equation} \label{conditdistribMx}
F_ {\mathcal M _x |E_0 } (z) = \frac {P [ \mathcal M _x \le z, {\cal X}(\tau _ \mu (x))=0]} {\pi ^ \mu  _0 (x)}= \frac {w_ {\mathcal M}(z)} {\pi ^ \mu  _0 (x)}.
\end{equation}

\begin{figure}
\centering
\includegraphics[height=0.4 \textheight]{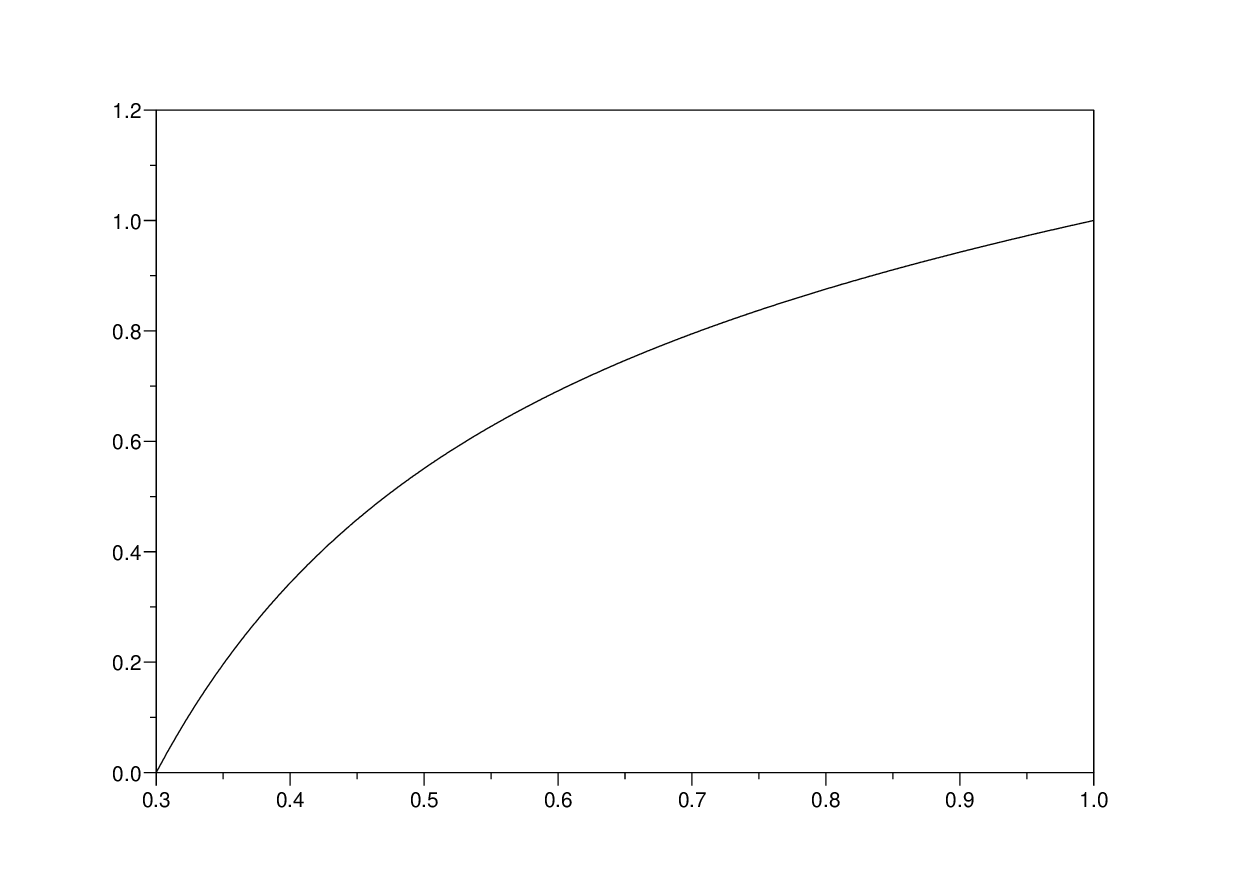}
\caption{An example of the graph of the conditional distribution of $\mathcal M _x ,$ given by \eqref{conditdistribMx}, as a function of $z \in (x,b),$  for
$b=1, \ \mu=1, \ r = 1, \ x_R= 1/4$ and $x=0.3 \ .$
}
\label{figmaximum}
\end{figure}

\begin{figure}
\centering
\includegraphics[height=0.4 \textheight]{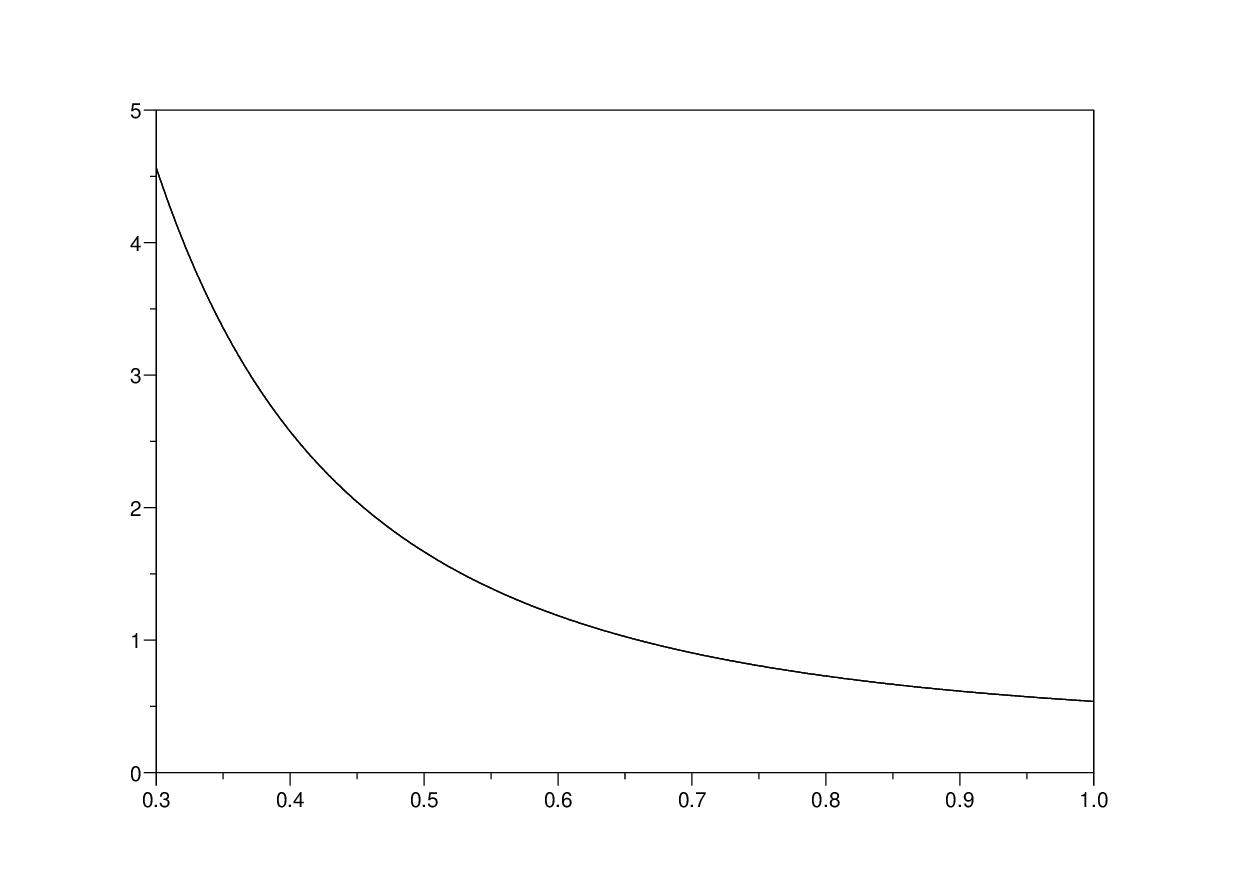}
\caption{An example of the graph of the conditional density of $\mathcal M _x $ given by \eqref{conditdensityMx}, as a function of $z \in (x,b),$  for
$b=1, \ \mu=1, \ r = 1, \ x_R= x=0.3 \ .$
}
\label{figmaximumdensity}
\end{figure}


\noindent For instance, if $x=x_R,$ \eqref{distrintM} becomes:
\begin{equation}
P \left [ \mathcal M _x \le z, {\cal X}(\tau _ \mu (x))  =0\right ]  =  \frac { \sinh \left ((z-x_R) \sqrt {\mu ^2 + 2r} \right ) } {\sinh \left ((z-x_R) \sqrt {\mu ^2 + 2r} \right ) + e^ {z \mu } \sinh \left (x_R  \sqrt {2r} \right )},
\end{equation}
and (see also Eq. (28) of \cite{guoyan:24}):
$$
\frac d {dz } P \left [ \mathcal M _x \le z, {\cal X}(\tau _ \mu (x))  =0\right ] =
$$
\begin{equation} \label{conditdensityMx}
= \begin{cases}
0 & \ \  z <x \\
\frac { e ^ {z \mu } \sinh \left (x \sqrt {\mu ^2 + 2r} \right ) \left [ \sqrt {\mu ^2 + 2r} \cosh \left ((z-x) \sqrt {\mu ^2 + 2r} \right ) - \mu \sinh \left ((z-x) \sqrt {\mu ^2 + 2r} \right ) \right ]}
{ \left [\sinh ((z-x) \sqrt {\mu ^2 + 2r}) + e^ {z \mu }\sinh (x  \sqrt {\mu ^2 + 2r})) \right ] ^ 2 } & \ \ x \le z \le b .
\end{cases}
\end{equation}
Thus, for $x = x_R$ the conditional density of the maximum till the exit time turns out to be:
\begin{equation}
f_{\mathcal M _x} (z) = \frac d {dz } P \left [ \mathcal M _x \le z | E_0 \right ] = \frac 1 {\pi _0 ^ \mu (x)} \frac d {dz } P \left [ \mathcal M _x \le z, {\cal X}(\tau _ \mu (x))  =0\right ], \ b \ge z \ge x ,
\end{equation}
where
$ \frac d {dz } P \left [ \mathcal M _x \le z, {\cal X}(\tau _ \mu (x))  =0\right ] $ is given by \eqref{conditdensityMx}.

In the Figure \ref{figmaximum} we report an example of the graph of the conditional distribution of $\mathcal M_x$  given by \eqref{conditdistribMx}, as a function of $z \in (x,b),$  for
$b=1, \ \mu=1, \ r = 1, \ x_R= 1/4$ and $x=0.3 \ .$
In the Figure \ref{figmaximumdensity} we report the corresponding graph of the conditional density of $\mathcal M_x$ given by \eqref{conditdensityMx}, as a function of $z \in (x,b),$
for $x=x_R =0.3$
(the other parameters are as in the Figure \ref{figmaximum}).
This conditional density behaves approximately as the exponential density, truncated to the interval
$(0.3, 1),$ namely $f_{ExpT}(z)= \theta e^{-\theta z}/(e^{-0.3 \theta } - e^{-\theta}), \ z \in (0.3,1),$ with $\theta = 4.5656;$ the mean turns out to be $0.49 \ .$
\bigskip

\noindent For drifted BM without resetting  $(r = 0)$ and $x=x_R ,$  we get
$$ \frac d {dz } P \left [\mathcal M _x \le z, {\cal X}(\tau _ \mu (x))  =0 \right ] = \frac {\mu  \sinh (\mu x) e^ {- \mu x}} {\sinh ^2 ( \mu z)}, \ z \ge x ,$$
(see e.g. Eq. (29) of \cite{guoyan:24}).
\bigskip

From \eqref{distrintM} with $\mu =0$ one gets: (cf. with Eq. (3.48) of \cite{abundo:FPA2023}, in the case $b= + \infty ):$
\begin{equation}
P[\mathcal M _x \le z, {\cal X}(\tau _ 0 (x))  =0] =  \frac { \sinh ((z-x_R) \sqrt {2r}) + \sinh ((x_R -x) \sqrt {2r}))} {\sinh ((z-x_R) \sqrt {2r}) + \sinh (x_R  \sqrt {2r}))}, \ z \ge x ,
\end{equation}
while it is zero for $ z <x .$  \par\noindent
If e.g. $x=x_R:$
\begin{equation}
P[\mathcal M _x \le z, {\cal X}(\tau _ 0 (x))  =0]  =  \frac { \sinh ((z-x_R) \sqrt {2r}) } {\sinh ((z-x_R) \sqrt {2r}) + \sinh (x_R  \sqrt {2r}))},
\end{equation}
and
\begin{equation}
\frac d {dz } P[\mathcal M _x \le z, {\cal X}(\tau _ 0 (x))  =0] =
\begin{cases}
0 & \ \  z <x \\
\frac { \sqrt {2r} \cosh ((z-x) \sqrt {2r} ) \sinh (x \sqrt {2r}) } {[\sinh ((z-x) \sqrt {2r}) + \sinh (x  \sqrt {2r}))]^2} & \ \ b \ge z \ge x .
\end{cases}
\end{equation}
Letting $r$ go to zero in the last equation, we get
\begin{equation}
\lim _{r \rightarrow 0} \frac d {dz } P \left [\mathcal M _x \le z, {\cal X}(\tau _ 0 (x))  =0 \right ] = \frac x {z^2}, \ z \ge x ,
\end{equation}
that is the density of the maximum displacement of BM without resetting (see e.g.
\cite{abundo:FPA2023}).
\bigskip

The minimum displacement of ${\cal X}(t)$ till the exit time, when starting from $x \in (0,b),$  is the  r.v.
${\rm m}_x = \min _ { t \in [0, \tau (x)]} {\cal X}(t)$
(of course we have $ 0 \le {\rm m}_x \le x).$ Also now, two cases are possible: in the first one ${\cal X}(t)$ first exits the interval $(0,b)$ through $0,$ and $ {\rm m}_x = 0,$ in the second case
${\cal X}(t)$ first exits the interval $(0,b)$ through $b$ and we can only say that ${\rm m}_x \le x.$ Thus,
 for $0 < z < x,$  $P[ {\rm m}_x >z , \ {\cal X}(\tau(x))= b]$ is nothing but the probability that  ${\cal X}(t)$ first exits the interval $(z,b)$ through $b.$
Therefore,
$w_{{\rm m}_x} (z):= P[{\rm m}_x >z , \ {\cal X}(\tau(x))= b]$ solves the differential boundary problem:
\begin{equation} \label{eqmindisplX}
\begin{cases}
 Lu(x) + r u(x_R) - ru(x) =0, \ x \in (z,b) \\
u(z) = 0, \ u(b) =1 ,
\end{cases}
\end{equation}
with $L u(x) = \frac 1 2 u''(x) + \mu u'(x).$
For $z=0$ its solution provides  $\pi_ b (x)$ again.
\par\noindent
Thus, the conditional survival  function of ${\rm m}_x$ is:
\begin{equation}
S_{{\rm m}_x | E_b} (z) = \frac 1 {\pi_ b (x)} P[  {\rm m}_x > z , \ {\cal X}(\tau(x))= b  ]= \frac {w_{{\rm m}_x} (x)}{\pi_b(x)},
\end{equation}
where $w_{{\rm m}_x} (x)$ is the solution of \eqref{eqmindisplX}. \par\noindent
By using Remark \ref{remark1}, we get for $0 < z < x:$
$$ w_{{\rm m}_x} (z)= P \left [{\rm m}_x  > z, {\cal X}(\tau _ \mu (x))  =b \right ]= $$
\begin{equation} \label{survintm}
  \frac { e^ {(b-z) \mu}\sinh \left ((x_R -z) \sqrt {\mu ^2 + 2r} \right ) - e^ {\mu (b-x)}\sinh \left ((x_R -x) \sqrt {\mu ^2 + 2r} \right  )} {\sinh \left ((b-x_R) \sqrt {\mu ^ 2 + 2r} \right ) + e^ {(b-z) \mu }\sinh \left ((x_R -z)  \sqrt {\mu ^2 + 2r} \right )} .
\end{equation}
Of course, $P \left [{\rm m}_x  >x, {\cal X}(\tau _ \mu (x))  =b \right ] =0$ and $P \left [{\rm m}_x  >0, {\cal X}(\tau _ \mu (x))  =b \right ] = \pi  ^ \mu _b (x)$ \par
\begin{figure}
\centering
\includegraphics[height=0.4 \textheight]{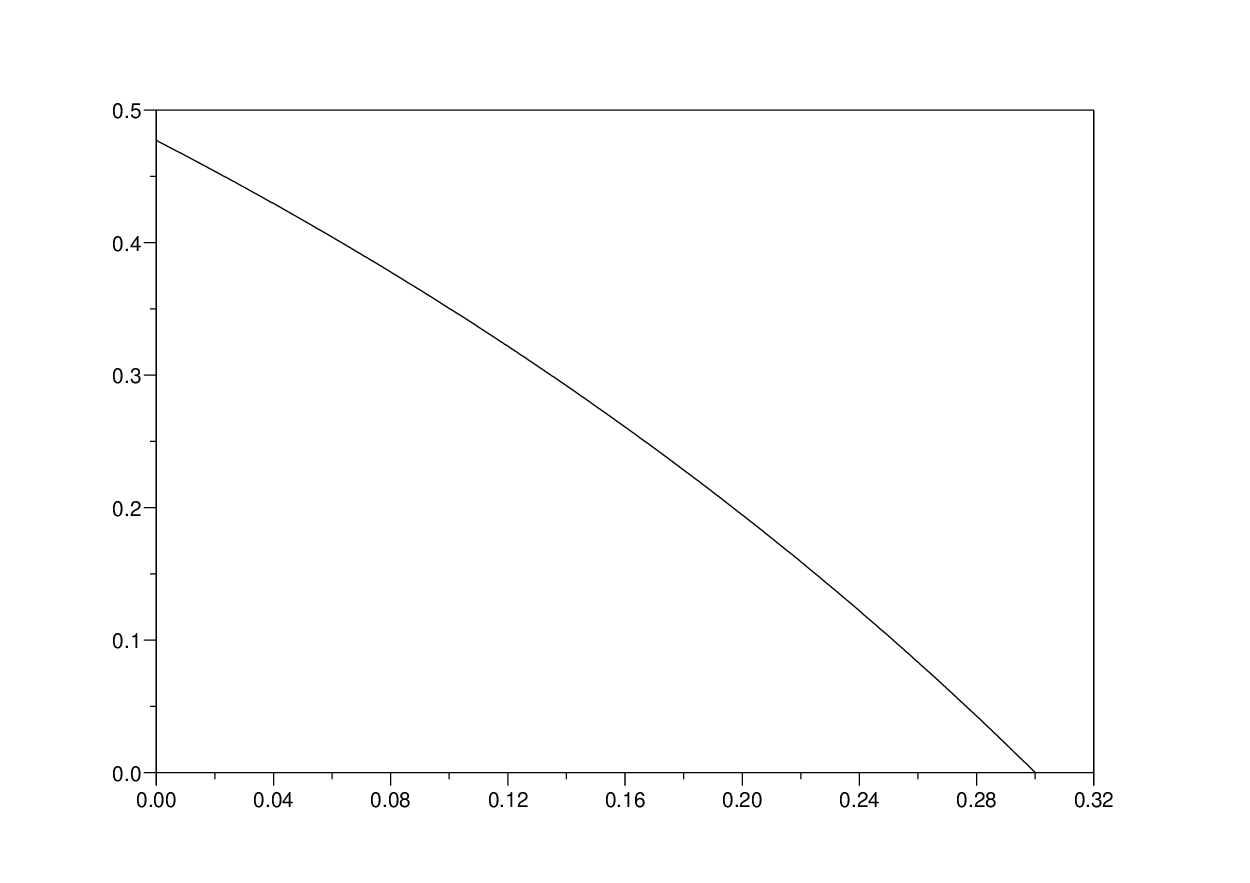}
\caption{Plot of $P \left [{\rm m}_x  > z, {\cal X}(\tau _ \mu (x))  =b \right ]$ as a function of $z \in (0,x),$  for
$b=1, \ \mu=1, \ r = 1, \ x_R= 1/4$ and $x=0.3 \ .$ It decreases from the value $\pi_b(x)=\pi _b(0.3)= 0.52286,$ attained at $z=0,$ to the value zero attained at  $z=x=0.3 \ .$
}
\label{figminprobab}
\end{figure}

Taking  $\mu =0$ in \eqref{survintm}, one gets (see also Eq. (12) of \cite{huang:24}):
\begin{equation}
w_{{\rm m}_x} (z) |_ {\mu =0 }=   \frac { \sinh \left ((z-x_R) \sqrt {2r} \right ) + \sinh \left ((x_R -x) \sqrt {2r} \right  )} {\sinh \left ((x_R -b) \sqrt {2r} \right ) + \sinh \left ((z-x_R )  \sqrt {2r} \right )}, \ 0 < z < x .
\end{equation}
In the Figure \ref{figminprobab}, we report an example of the graph of $P \left [{\rm m}_x  > z, {\cal X}(\tau _ \mu (x))  =b \right ]$ for
$b=1, \ \mu=1, \ r = 1, \ x_R= 1/4$ and $x=0.3 \ .$ It decreases from the value $\pi_b(x)=\pi _b(0.3)= 0.52286,$ attained at $z=0,$ to the value zero attained at  $z=x=0.3$


\section{The inverse first-exit time (IFET) problem for BM with resetting}
In this section, we study the IFET problem  for (undrifted) BM with resetting, ${\cal X}(t),$ in a bounded interval $(0,b).$ Therefore, the underlying diffusion is
$X(t) = \eta + W_t ,$
where the initial position $\eta = {\cal X}(0) $ is supposed to be randomly distributed in $(0,b)$ and independent of  ${\cal X}(t),$ whereas the reset rate $r >0$ and the reset position $x_R \in (0,b)$ are fixed.
Actually, we omit to treat the case of drifted BM with resetting $(\mu, \ r  \neq 0),$ because of  heavier calculations.\par
We recall the terms of the IFET problem, mentioned in the Introduction. \par\noindent
Let $\tau(x)= \tau _0(x)$ be the FET of $\mathcal X(t)$ from $(0,b),$ under the condition that $\eta  =x \in (0,b),$ namely
\begin{equation}
\tau (x) = \inf \{ t>0: \mathcal X(t) \notin (0,b)  | \eta = x \},
\end{equation}
and let
$\tau$ the (unconditional) FET of $\mathcal X(t)$  i.e.
\begin{equation}
\tau = \inf \{ t>0: \mathcal X(t) \notin (0,b) \}.
\end{equation}
For a given  distribution function $F (t)$ on the positive real axis, or equivalently for a given FPT density $f(t) = F'(t), \ t >0,$ the IFET problem
consists in finding the density $g$ of the
random initial position $\eta \in (0, b),$  if it exists, such that $P[\tau \le t ] = F(t), \ t >0.$ The function $g$ is called a solution to the IFET problem
for $\mathcal X(t);$  in fact, the uniqueness of the solution is not guaranteed (see Example 5).
For BM without resetting $(r=0),$ the IFET problem was studied in \cite{abundo:stapro13}.\par\noindent
If $g(x)$ is the density of the random initial position $\eta ,$
by using the explicit form \eqref{LTtauBMreset} of the LT $M_0(x, \lambda)$ of $\tau (x),$
we obtain
the LT of the FET $\tau:$
$$ E \left [e^ {- \lambda \tau} \right ] = \int _0 ^b g(x) M_0(x, \lambda) dx $$
$$=  \int _0 ^b g(x) \left \{ 1 - \frac  {\lambda  \left [ \sinh ( \alpha _ \lambda  b) - \sinh (\alpha _ \lambda x ) - \sinh (\alpha _ \lambda (b-x) ) \right ] }    {{\cal A}_ \lambda } \right \} dx $$
\begin{equation} \label{explLTtauBMreset}
 = 1 - \frac \lambda {{\cal A}_ \lambda } \int _0 ^b g(x)  \left [ \sinh (\alpha _ \lambda b ) - \sinh (\alpha _ \lambda x ) - \sinh (\alpha _ \lambda (b-x) ) \right ]   dx ,
\end{equation}
or also
\begin{equation} \label{explLTtauBMresetprimo}
E \left [e^ {- \lambda \tau} \right ] = 1 -  \frac 1 {{\cal A}_ \lambda } \Big \{ \lambda \sinh (\alpha _ \lambda b) - \lambda E \left [ \sinh (\alpha _ \lambda \eta ) + \sinh (\alpha _ \lambda (b-\eta ) ) \right ] \Big \} ,
\end{equation}
where
$$
{\cal A}_ \lambda  =  \lambda  \sinh (\alpha _ \lambda b ) + r \sinh (\alpha _ \lambda x_R ) + r \sinh (\alpha _ \lambda (b-x_R) ) ,$$
and
$$
 \alpha _ \lambda = \sqrt {2( \lambda +r)}
$$
(for simplicity, we omit to explicit the dependence of ${\cal A}_ \lambda $ on $b, r, x_R $ and of $\alpha _ \lambda $ on $r).$
\par\noindent
Let us denote by $\widehat g (\theta)= \int _0 ^b e^{- \theta x} g(x) dx $  the (possibly bilateral) LT of $g(x).$ Then, from \eqref{explLTtauBMreset}  we obtain the following:
\begin{Proposition} \label{Prop2}
Let $\mathcal X(t)$ be  (undrifted) BM with resetting, starting from the random initial position $\eta \in (0, b),$ which is supposed to be independent of $\mathcal X(t),$ and let $f(t), \ t >0,$ be a given  FPT density. Then, if there
exists a  solution $g$  to the IFET problem  for $\mathcal X(t),$ the following equality holds:
\begin{equation} \label{LTtau}
\widehat f  (\lambda ) = 1-  \frac \lambda {{\cal A}_ \lambda} \left \{ \sinh (\alpha _ \lambda b) - \frac 1 2 \left [ \left (e^{\alpha _ \lambda b} -1 \right ) \widehat g (\alpha _ \lambda)
+ \left (1-e^{- \alpha _ \lambda b} \right ) \widehat g (- \alpha _ \lambda ) \right ]   \right \}, \ \lambda >0,
\end{equation}
 $(\alpha _ \lambda = \sqrt {2(\lambda +r)}).$ \par
If one requires that $g$ is symmetric with respect to $x=b/2,$ that is, $g \left ( \frac b 2  -u \right )= g \left ( \frac b 2  +u \right ), \ u  \in (0, b/2),$ one finds $\widehat g (\theta)= e^ {-b \theta} \widehat g (- \theta), \ \theta >0,$ and so \eqref{LTtau} becomes:
\begin{equation} \label{LTf}
\widehat f  (\lambda ) =1 - \frac \lambda {{\cal A}_ \lambda } \left [ \sinh(\alpha _ \lambda b )-(e^{\alpha _ \lambda} -1) \widehat g (\alpha _ \lambda) \right ],
\end{equation}
or
\begin{equation} \label{LTg}
\widehat g  (\theta ) =  \frac {\Big [ \left (\frac {\theta ^2} 2 -r \right ) \sinh( \theta b)+r s_ \lambda \Big ]  \widehat f \left (\frac {\theta ^2} 2 -r \right )
- r s_ \lambda} {(\frac {\theta ^2} 2 -r) (e^ {\theta b} -1) },
\end{equation}
where
\begin{equation}
s_ \lambda  = \sinh(\alpha _ \lambda x_R)+ \sinh(\alpha _ \lambda (b-x_R) )
\end{equation}
(we omit the dependence of $s_ \lambda $ on $b, r, x_R ).$
\hfill  $\Box$
\end{Proposition}
\begin{Remark} \label{oss}
For $r=0$ Eq. \eqref{LTg} coincides with Eq. (2.2) of \cite{abundo:stapro13}. \par\noindent
If we set $\bar T = \max _ {x \in (0,b)} E[\tau (x)],$ then from
$E[\tau ] = \int _0 ^b E[\tau  (x)] g(x) dx, $ we get the compatibility condition
\begin{equation}
E[\tau ] \le \bar T ,
\end{equation}
which is necessary so that a solution to the IFET problem exists. \par
Formula
 \eqref{LTtau} provides a functional relation between the LT of $f  $ and that of $g.$
Once $\widehat g$ has been found, such that it satisfies  \eqref{LTtau}, it may be that $\widehat g$ is not the LT of the probability density function of a random variable. In this case, a solution to the IFET problem does not exist. This is the reason way Proposition \ref{Prop2} is formulated in a conditional form.
 \par\noindent
In the case when $g$ is symmetric with respect to $b/2,$ \eqref{LTg} allows to write
$\widehat g  (\theta )$ in terms of a function of $\widehat f (\theta )$ and other parameters.
If $\widehat f$ is an analytic function, then
$\widehat g$ is also analytic in the interval $(- \sqrt {2r}, \sqrt {2r}), $ and so \eqref{LTg}
uniquely identifies the density $g$, and hence the distribution of $\eta;$ therefore, if there is a solution to the IFET problem, then it is unique. \par\noindent
Thus, we conclude that, if $\widehat f $ is analytic, then the solution $g$ to the IFET problem is unique, under the constraint that it is sought in the class of densities  that are symmetric with respect to $b/2.$
\end{Remark}
\bigskip

Now, we  will prove the existence of a solution $g$ to the IFPT problem
for a class of FPT densities $f.$
For the sake of simplicity, we limit ourselves to the case when $b=1, \ r=1, \ x_R = 1/2 $  and $g$ is  symmetric with respect to $b/2 = 1/2.$
\par\noindent
For any integer $k \ge 0,$ set $J_k(\theta) = \int _{0} ^1 e^ { - \theta x } x^k dx;$ as easily seen,
$J _0 (\theta )= \frac {1-e^ {- \theta}} \theta $ and the recursive relation
$ J _ {k} (\theta ) = \frac k \theta J_ {k-1} (\theta) - \frac {e^ {-\theta}} \theta, \ k \ge 1$ allows to
calculate $J _k (\theta),$ for every $k.$ \par
The following Proposition gives a sufficient condition, in order that there exists a
solution to the  IFPT problem for BM with resetting.

\begin{Proposition} \label{existenceproposition}
Let ${\cal X}(t)$ be (undrifted) BM with resetting, and let $b=1, \ r=1, \ x_R = 1/2 .$
Suppose that the LT  of $f(t)$ has the form:
\begin{equation} \label{laplacedensityclass}
\widehat f ( \lambda )= \sum _ {k=1} ^ N p_k \widehat f _ {k}( \lambda), \ \lambda >0,
\end{equation}
where:
\begin{equation} \label{ESLTfk}
\widehat f _ {k} ( \lambda ) = \frac {2 \sinh(\alpha _ \lambda /2) + \lambda (e^ {\alpha _ \lambda } -1 ) \widehat g _ {k}(\alpha _ \lambda) } { \lambda \sinh (\alpha _ \lambda) + 2 \sinh (\alpha _ \lambda /2 ) } \ \
 (\alpha _ \lambda = \sqrt {2( \lambda +r)}),
\end{equation}
the numbers $p_k \ge 0 $ are chosen in such a way that $\sum _{k=1} ^N p_k =1,$ \par\noindent
and
\begin{equation} \label{LTgkProp}
\widehat g_{k}(\lambda)= \frac { (2k+1)!} {(k!)^2} \sum _ {i=0} ^k (-1) ^i
\begin{pmatrix}
  k \\ i
\end{pmatrix}
 J _ {k+i} (\lambda), \ k \ge 0.
\end{equation}
Then,
there exists a solution $g$ to the IFPT problem for ${\cal X}(t)$,  corresponding to
the FPT density $f,$ and
it results:
\begin{equation} \label{densityclass}
g(x)= \sum _ {k=1} ^N p_kg_{k} (x).
\end{equation}
Note that the density $g_k$ corresponding to the LT \eqref{LTgkProp} is the Beta density in the interval $(0,1)$ with parameters $\alpha = \beta = k+1,$ i.e.
\begin{equation}
g_k(x)= \frac { (2k+1)!} {(k!)^2} x^k (1-x) ^ k, \ k \ge 0 .
\end{equation}
\end{Proposition}
{\it Proof.}
First, we observe that,
if the density $g_{k}, \ k=1, \dots N ,$ is a solution to the IFPT
problem corresponding to the FPT
density $f_{k},$  then $g(x)= \sum _ {k=1}^N p_k g_{k}(x) ,$ with $p_k
\ge 0$ and $\sum _{k=1}^ N p_k =1,$ solves the IFPT problem for
the FPT density $f(t) = \sum _k p_k f_{k} (t).$ Thus, it is enough to
verify that $g_{k}$ is a solution to the IFPT problem corresponding to
$  f _{k}  .$ Since a simple calculation shows that the LT of $g_k$ is:
\begin{equation}
\widehat g_{k}(\lambda)= \frac { (2k+1)!} {(k!)^2} \sum _ {i=0} ^k (-1) ^i
\begin{pmatrix}
  k \\ i
\end{pmatrix}
 J _ {k+i} (\lambda), \ \lambda >0,
\end{equation}
the verification follows by inserting $\widehat g _ {k}$ and $\widehat f _ {k}$ given by \eqref{ESLTfk} into \eqref{LTf}, with $b=1, \ r=1, \ x_R = 1/2 ,$ because $g_k(x)$ is  symmetric with respect to $x= 1/2 .$

\par \hfill  $\Box$

In the next subsection, we show some explicit examples of solutions to the IFET problem for (undrifted) BM with resetting.
\subsection{Some examples of solutions to the IFET problem}
{\bf Example 1.} Let ${\cal X}(t)$ be BM with resetting, and let $b=1, \ r=1, \ x_R = 1/2;$ suppose that the LT of the FET $\tau$ is:
\begin{equation} \label{LTEx1}
\widehat f (\lambda) =
2 \  \frac {\sqrt {2(\lambda +1)} \sinh \left (\frac 1 2 \sqrt {2(\lambda +1)} \right ) + \lambda \left ( \cosh \left (\sqrt {2(\lambda +1)} \right )- 1 \right )} { \sqrt {2(\lambda +1)} \left [ \lambda \sinh \left ( \sqrt {2(\lambda +1)} \right ) + 2 \sinh \left (\frac 1 2 \sqrt {2(\lambda +1)} \right )  \right ] }, \ \lambda >0 .
\end{equation}
We search for a solution to the IFET problem for $\mathcal X (t)$ in the set of probability densities $g$ in $(0,1)$ which are symmetric with respect to $1/2;$ we find that a solution $g$  is the uniform density in the interval $(0,1),$ i.e. $g(x)= \mathbb{I}_ {(0,1)} (x).$
Since $\widehat f (\lambda)$ is analytic, from Remark \ref{oss} it follows that this is the only solution in the class of densities which are symmetric with respect to $1/2.$  \par\noindent
To verify this, we can  use \eqref{LTtau} or \eqref{LTg}: it is sufficient to substitute the various quantities
 and use that for the uniform density $g$ in $(0,1)$ one has
  $\widehat g( s)= (1-e^{-s} ) / s.$ \par\noindent
The LT \eqref{LTEx1} cannot be inverted in closed form to obtain the corresponding FET density $f;$ however, we can get some qualitative characteristics of the distribution having density $f.$ In fact, from $\widehat f  (\lambda )$ we easily obtain all the moments $m_k$ and the central ones $\mu_k $ of the corresponding distribution.
By rounding the values to the third decimal digit, the mean of the distribution turns out to be $E[\tau]= 0.175, $ while
$ E[\tau ^2]= 0.074, \ \ E[\tau ^3]= 0.047, \ \ E[\tau ^4]= 0.041 \ .$
Thus,
the first three central moments are
$  \mu_2 = 0.043, \ \mu _3 = 0.019, \ \mu _4= 0.018 , $
from which skewness  $\gamma _1 = \frac {\mu _3 } {\mu_2 ^{3/2}}  = 2.139,$ and excess kurtosis coefficient $\gamma _2 = \frac {\mu _4} { \mu _2 ^2 } -3=
6.632$ follow. Since skewness $\gamma _1 $ is positive,
the tail of the distribution is on the right side;
moreover, the density $f(t)$ tends to zero, as $ t \rightarrow + \infty,$ more slowly than the normal density does,  because excess kurtosis $\gamma _2 $ is positive. \par
\begin{figure}
\centering
\includegraphics[height=0.28 \textheight]{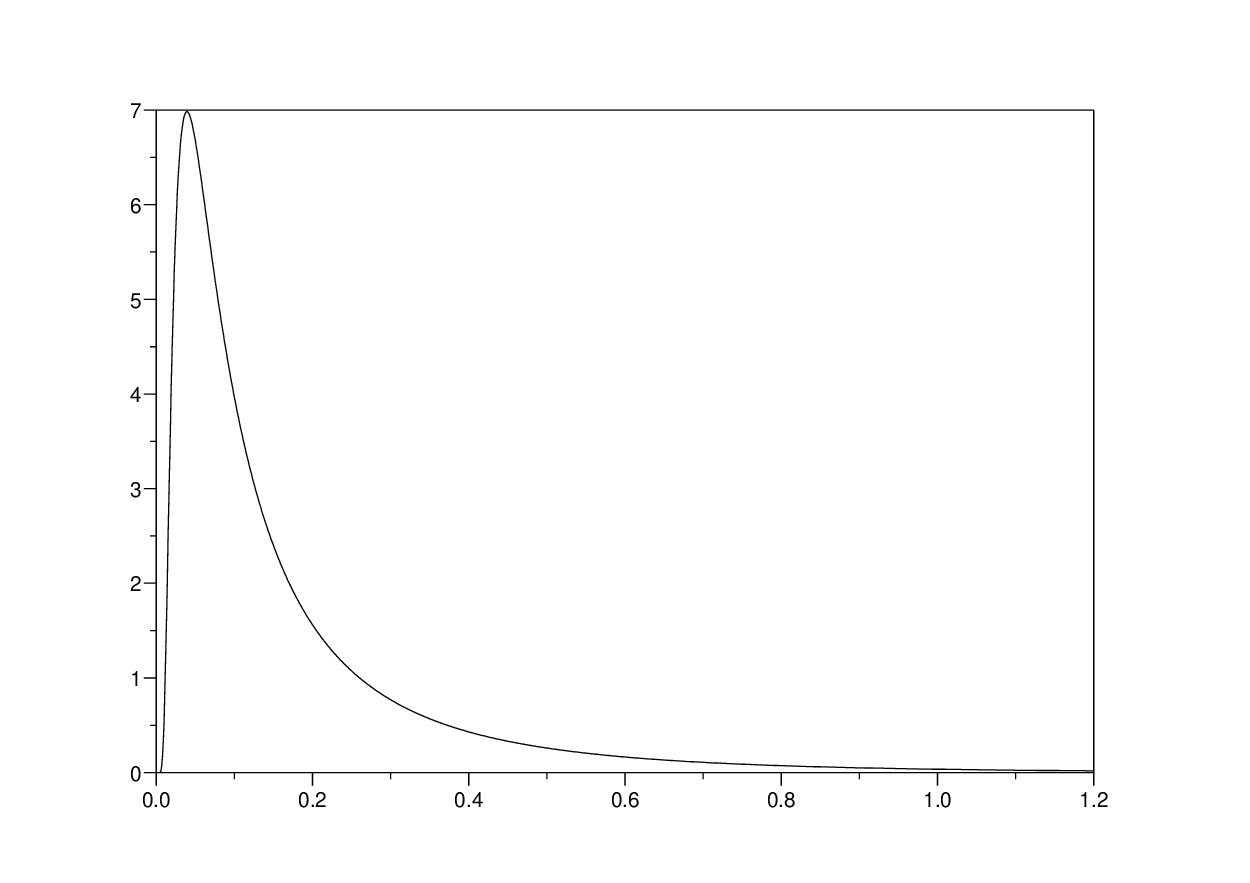}
\caption{Qualitative graph of the FET density  $f(t)$ corresponding to \eqref{LTEx1}
}
\label{InvGauss}
\end{figure}

\noindent
A qualitative graph of the density $f(t),$ as a function of $t \ge 0,$ is shown in the Fig. \ref{InvGauss}. \bigskip

\bigskip

\noindent {\bf Example 2.}
Let ${\cal X}(t)$ be BM with resetting, and  let $b=1, \ r=1, \ x_R = 1/2,$ and suppose that the LT of the FET $\tau$ is:
\begin{equation}
\widehat f  ( \lambda ) = \frac {2 \sinh \left ( \frac 1 2 \sqrt {2(\lambda +1)} \right ) + \lambda (e^ {\sqrt {2(\lambda +1)}} -1 ) \widehat g _ {1}(\sqrt {2(\lambda +1)}) } { \lambda \sinh \left (\sqrt {2(\lambda +1)} \right ) + 2 \sinh \left (\frac 1 2 \sqrt {2(\lambda +1)} \right ) }, \ \lambda >0.
\end{equation}
Then, a solution $g$ to the IFET problem for ${\cal X}(t)$ is the Beta density  $g_1(x)= 6 x (1-x) \mathbb{I}_{(0,1)}(x).$
This is a special case of Proposition \ref{existenceproposition} for $N=k=1.$
\bigskip

\noindent {\bf Example 3.}
Let ${\cal X}(t)$ be BM with resetting, and  let $b=1, \ r=1, \ x_R = 1/2,$ and $\alpha, \beta >0;$ suppose that the LT of the FET $\tau$ is:
$$
\widehat f (\lambda) =
1 - \left [ \frac \lambda {\lambda \sinh(\alpha _ \lambda) + 2 \sinh(\alpha _ \lambda /2)} \right ]  \times $$
\begin{equation} \label{LTEx2}
\times \left \{  \sum _ {k=0}^ \infty \frac {\alpha _ \lambda ^ {2k+1}} { (2k+1)!} (\cosh(\alpha _ \lambda) -1) \frac { B(\alpha +k, \beta )} {B( \alpha, \beta)} -
\sum _ {k=1}^ \infty \frac {\alpha _ \lambda ^ {k}} { (2k)!} \sinh (\alpha _ \lambda)\frac { B(\alpha +k, \beta )} {B( \alpha, \beta)} \right \} ,
\end{equation}
where $\alpha _ \lambda = \sqrt {2( \lambda +r)}$ and
$$B( \alpha, \beta) = \frac {\Gamma (\alpha) \Gamma (\beta)} {\Gamma( \alpha + \beta)}.$$
Then, a solution $g$ to the IFET problem for ${\cal X}(t)$ is the Beta density in $(0,1)$ with parameters $\alpha$ and $\beta,$  namely
$g(x) = \frac { \Gamma (\alpha + \beta)} {\Gamma(\alpha) \Gamma ( \beta)} x^{\alpha -1} (1-x) ^{\beta -1} \mathbb{I}_{(0,1)}(x).$
Example 1 and Example 2 are special cases, when $\alpha = \beta =1,$ and $\alpha = \beta =2,$ respectively.
\par\noindent
To verify this by means of Proposition \ref{Prop2} it is sufficient to substitute the various quantities into \eqref{explLTtauBMresetprimo} or \eqref{LTtau}; it is convenient  to use  that, if $\eta$ has Beta density in $(0,1),$ one has (see e.g. \cite{gupta}):
\begin{equation} \label{BetaMGF}
 E \left [e^ {s \eta} \right ]= \sum _{k=0} ^\infty \frac {s^k} {k!} \frac { B(\alpha +k, \beta )} {B( \alpha, \beta)}.
\end{equation}
Thus:
\begin{equation}
E \left [\sinh (s \eta) \right ] = \sum _ {k=0} ^ \infty \frac {s ^ {2k+1}} { (2k+1)!} \frac { B(\alpha +k, \beta )} {B( \alpha, \beta)},
\end{equation}
and
\begin{equation}
E \left [\sinh (s (b - \eta )) \right ] =   \sum _ {k=0} ^ \infty \frac {s ^ {2k}} { (2k)!} \sinh ( s b)\frac { B(\alpha +k, \beta )} {B( \alpha, \beta)}
- \sum _ {k=0} ^ \infty \frac {s ^ {2k+1}} { (2k+1)!} \cosh ( s b) \frac { B(\alpha +k, \beta )} {B( \alpha, \beta)}.
\end{equation}

\bigskip

\noindent {\bf Example 4.} Let ${\cal X}(t)$ be BM with resetting, and let $b=1,  \ x_R = 1/2,$ and $ \gamma $ a fixed positive number; suppose that the LT of the FET $\tau$ is, for $\lambda \neq \gamma ^2 /2 -r:$
$$
\widehat f (\lambda) =
\frac 1 {\lambda \sinh (  \sqrt {2(\lambda +r)}) +r \left [\sinh ( x_R  \sqrt {2(\lambda +r)})+ \sinh((1-x_R)\sqrt {2(\lambda +r)}) \right ]}  \times
$$
\begin{equation}
\Biggl \{ r \left [\sinh ( x_R  \sqrt {2(\lambda +r)})+ \sinh((1-x_R)\sqrt {2(\lambda +r)}) \right ] +
\end{equation}
$$
+ \frac  {\lambda \gamma } {2 (e^ \gamma -1)} \left [ \frac { \left (e^  {\sqrt {2(\lambda +r)}} -1  \right ) \left (e^ \gamma - e^ {- \sqrt {2(\lambda +r)}} \right )    }  {\gamma + \sqrt {2(\lambda +r)} }
+ \frac { \left (1- e^ {- \sqrt {2(\lambda +r)}}  \right ) \left (e^ \gamma - e^ {\sqrt {2(\lambda +r)}} \right ) }{\gamma - \sqrt {2(\lambda +r)}}
\right ] \Biggr \}. $$
Then, a solution $g$ to the IFET problem for ${\cal X}(t)$  is the exponential density with parameter $\gamma >0,$  truncated to the interval $(0,1),$ i.e.
$g(x)= \frac {\gamma e^ {- \gamma x}} {1- e^ {- \gamma}}\mathbb{I}_ {(0,1)} (x).$
To verify this by means of Proposition \ref{Prop2} it is sufficient to substitute the various quantities into \eqref{LTtau}, and to use that for this density $g$ one has
  $\widehat g( s)= \frac {\gamma ( e^ \gamma - e ^ {-s})} {(e^ \gamma -1 )( \gamma + s)}.$ \par\noindent

\bigskip

\noindent {\bf Example 5.} Let ${\cal X}(t)$ be BM with resetting, and let  $b=1, \ r=1, \ x_R = 1/2 ;$ suppose that the LT of the FET $\tau$ is given again by
\eqref{LTEx1} of Example 1. Now, we seek for a solution $g$ to the IFET problem, without the constraint on the symmetry of $g$ with respect to $1/2;$ therefore, we must use
\eqref{LTtau}. Then, we obtain that a solution to the IFET problem for ${\cal X}(t)$  is the density $g(x)= 2x \mathbb{I}_ {(0,1)} (x).$
This is easily verified by substituting  the various quantities into \eqref{LTtau}, and using that for this density $g$ one has
$\widehat g( s) = \frac 2 {s^2} (1- e^s (1-s)).$ Thus, for the LT $\widehat f$ given by \eqref{LTEx1}, we have found two different solutions to the IFET problem:
this last solution $g$ and that of Example 1. \par\noindent
This example shows that there is not uniqueness of the solution $g$ to the IFET problem, unless one does not introduce constraints on the class of  densities $g.$

\bigskip

\noindent In the next example the solution to the IFET problem is a discrete density.
\bigskip

\noindent {\bf Example 6.} Let ${\cal X}(t)$ be BM with resetting, and let  $b=2,\  r=1, \ x_R = 1/2 ;$ we suppose that the random initial  position $\eta$ takes integer values in the set $  \{ 0,\ 1,\ 2 \},$
and that the LT of the FET $\tau$ is:
\begin{equation}
\widehat f (\lambda) =
2 \  \frac { \sinh \left (\frac 1 2 \sqrt{ 2 (\lambda +1) } \right ) +\frac \lambda 3 \left [ \sinh \left (\sqrt{ 2 (\lambda +1) } \right ) + \sinh \left (2 \sqrt{ 2 (\lambda +1) }\right ) \right ]} { \lambda \sinh \left (2 \sqrt{ 2 (\lambda +1) } \right ) + 2 \sinh \left (\frac 1 2  \sqrt{ 2 ( \lambda +1) } \right )  }.
\end{equation}
This LT cannot be inverted in closed form to obtain the corresponding FET density $f;$ however,
the mean of the distribution turns out to be $E[\tau]= 0.991, $ and
$ E[\tau ^2]= 8.148 \ . $
We seek for a solution $g$ to the IFET problem, with the constraint that $g$ is symmetric  with respect to $x=1;$ therefore, we can use
\eqref{LTg}. By substituting  the various quantities into \eqref{LTg}, we find that
a solution $g$ to the IFET problem for ${\cal X}(t)$ with this constraint,  is the discrete uniform density in the set $\{ 0, 1, 2 \},$ namely $g(x) = \frac 1 3, \ x \in \{0,\ 1,\ 2 \},$ being
$\widehat g ( \lambda )= \frac 1 3 (1+ e^ {- \lambda} + e^ {-2 \lambda}  ).$

\section{Conclusions and Final Remarks}
We studied the direct and  inverse first-exit time (FET) problems for a one-dimensional diffusion process $\mathcal X(t)$ with resetting, obtained from an
underlying temporally homogeneous diffusion  $X(t)$
driven by the SDE
$ dX(t)=  \mu(X(t)) dt + \sigma (X(t)) d {W_t} ,$
where $W_t$ is a standard Brownian motion, and the drift $\mu (\cdot)$ and diffusion coefficient $\sigma (\cdot)$ are regular enough functions, such that there exists a unique strong solution of the SDE.
The process $\mathcal X(t)$ starts from $\mathcal X(0)= \eta \in (0,b) $ and it is subject to reset to the position $x_R$ according to a homogeneous Poisson process with rate $r >0.$
\par\noindent
Thus, for any $C^2$ function $u(x),$ the infinitesimal generator $ \mathcal L$ of ${\cal X}(t)$ is given by
$
{\cal L}u(x) = \frac 1 2 \sigma ^2(x) u''(x) + \mu (x) u'(x) +r (u(x_R) -u(x)) = L u(x) +r (u(x_R) -u(x)) ,
$
where $Lu(x)= \frac 1 2 \sigma ^2(x) u''(x) + \mu (x) u'(x)$ represents the  ``diffusion part'', i.e. that concerning the diffusion  $X(t).$ \par
As regards the direct FET problem, for  fixed $\eta = x \in (0,b)$ we investigated the statistical properties of the FET $\tau (x) $ of $\mathcal X(t)$ from the interval $(0,b),$ where $\tau (x) = \min \{t>0: \mathcal X(t) \notin (0,b) | \mathcal X(0) =x \};$ moreover, we studied the probability distribution of the first exit area FEA $A (x),$ namely the area swept out by ${\cal X}(t)$  till the time $\tau (x),$ and the probability distributions of the maximum and minimum displacement of ${\cal X}(t)$ for  $t \in (0, \tau (x)).$ \par\noindent
In particular, we stated ODEs with boundary conditions for the exit probability of  $\mathcal X(t)$ from the left and right ends of the interval $(0,b),$ for the Laplace transforms of $\tau (x)$ and $A(x),$ for the single and joint moments of $\tau(x)$ and $A(x),$ and for
the probability distributions of the maximum and minimum displacement of ${\cal X}(t)$ till the FET, all in terms of the infinitesimal generator ${\cal L}$ of $\mathcal X (t).$
\par
In the case of drifted Brownian motion with resetting, namely when the diffusion and drift coefficients are $\sigma (x) =1$ and
$\mu (x) = \mu ,$
the generator turns out to be  $\mathcal L u(x)= \frac 1 2  u''(x) + \mu u'(x) +r (u(x_R) -u(x)) ,$ and
the ODEs were explicitly solved in terms of elementary functions, except that concerning the Laplace transform of $A(x),$ that admits a solution only in terms of special functions. \par
As for the inverse first-exit time (IFET) problem for (undrifted) Brownian motion with resetting, we supposed that the starting position $\eta$ was randomly distributed in $(0,b),$ and we denoted by $\tau$ the FET of $\mathcal X(t)$ from the interval $(0,b);$ then, for a given distribution function $F(t)$ on the the time $t$ axis, the IFET problem consisted in finding the probability density $g$ of  $\eta,$  if existing, such that
the FET of $\mathcal X(t)$ from $(0,b)$ has distribution $F(t),$ namely
$P[\tau \le t ] = F(t), \ t >0.$ The density $g$ was called a solution to the  IFET problem for $\mathcal X(t);$
it turned out to be unique, under certain conditions.
We reported several explicit examples of solutions to the IFET problem, for Brownian motion with resetting. \par\noindent
The case when $b= + \infty$ was studied in \cite{abundo:FPA2023}, while the two-barrier IFET problem for diffusions without resetting was studied in \cite{abundo:stapro13}. \par
\par
A different type of IFET problem for diffusions with resetting, consisting in finding  the density $g$ of the starting position $\eta \in (0,b)$ corresponding to an assigned
mean value of the FET $\tau,$ was studied in \cite{abundo:TPMS2024}; of course, in that case the solution $g$ was not unique. \par
We remark that the inverse problem here considered concerns randomization in the starting point $\eta $ of $\mathcal X(t).$
More generally, one could introduce randomization in the reset position $x_R$ (taking fixed $\eta),$ or both in the starting point $\eta$ and in the reset position $x_R,$ and then
study the corresponding IFET problems, where now a solution is the joint density of $(\eta, \ x_R),$ if it exist.
\par
The feature of the present paper was to characterize several quantities concerning direct and inverse FET problems, as solutions of ODEs with boundary conditions, which can be explicitly
solved in terms of elementary functions. \par\noindent
Although the calculations were developed for drifted Brownian motion with resetting,
in principle they can be carried on for any one-dimensional diffusion with resetting, obtained from $X(t)$ driven by the SDE \eqref{diffusion};
it suffices to substitute the corresponding generator $ Lu(x)= \frac 1 2 \sigma ^2(x) u''(x) + \mu (x) u'(x)$ in all the ODEs.
\par\noindent
For instance,  one can study the case when  the underlying diffusion $X(t)$ is {\it conjugated} to Brownian motion, namely
there exists an increasing differentiable function $v(x)$ with $v(0) = 0,$ such
that $X(t)= v^{-1} \left (W_t + v( \eta ) \right ),$ for any $t \ge 0.$ (see \cite{abundo:stapro12}).
Actually, if  $\mathcal X (t)$ is obtained from a diffusion $X(t)$ which is conjugated to Brownian motion via the function $v$, the direct and inverse FET problems of $\mathcal X (t)$ are easily reduced to those of BM with resetting in the interval $(0, v(b)),$  starting from $v(\eta ).$ \par\noindent
Another diffusion with resetting that can be reduced to drifted BM with resetting is the process $\mathcal X (t),$ whose underlying diffusion is Geometric Brownian motion, that is
driven by the SDE $dX(t)= \nu X(t) dt + \sigma X(t) dW_t, \ \nu, \sigma >0.$  For other examples, see e.g. \cite{abundo:TPMS2024}. \par
Moreover, the arguments of this article can also be applied, for example, to the Ornstein-Uhlenbeck process
with stochastic resetting (see e.g. \cite{dubey}), namely when the underlying diffusion $X(t)$  is driven by the SDE
$dX(t) = - \nu X(t) + \sigma dB_t,$ for positive constants $\nu$ and $\sigma.$ Of course, in this case the corresponding differential equations to obtain the various quantities are rather complicated, however they can be solved in principle.
\par
Our study was motivated by the fact that, as in the case without resetting, direct and inverse problems for the FET of a diffusion process are  worthy of attention, because they have notable applications in several applied fields, for instance
in biological modeling concerning neuronal activity,  queuing theory, and  mathematical finance.

\bigskip

\noindent {\bf Acknowledgments:}
The author belongs to GNAMPA, the Italian National Research Group of INdAM; he also acknowledges
the MUR Excellence Department Project MatMod@TOV awarded to the Department of Mathematics, University of Rome Tor Vergata, CUP E83C23000330006 .


\begin{thebibliography}{99}
\bibitem {Abramowitz}
Abramowitz, M., Stegun, I.A., 1965. \newblock Handbook of mathematical functions: With formulas, graphs, and mathematical tables. Dover, New York

\bibitem {abundo:TPMS2024}
Abundo, M., 2024. \newblock  Inverse first-passage problems of a diffusion with resetting.
\newblock {Theor. Probability and Math. Statist }. To appear.

\bibitem {abundo:FPA2023}
Abundo, M., 2023. \newblock  The first-passage area of a Wiener process with stochastic resetting.
\newblock {Methodol Comput Appl Probab } 25:92 https://doi.org/10.1007/s11009-023-10069-4

\bibitem {abundo:OUarea}
Abundo, M.,  2023. \newblock The first-passage area of Ornstein-Uhlenbeck process revisited.
\newblock Stochastic Analysis and Applications 41(2): 358--376. https://doi.org/10.1080/07362994.2021.2018335.
\bigskip

\bibitem {abundo:CC22}
Abundo, M., 2022. \newblock  Some examples of solutions to an inverse problem for the first-passage place of a jump-diffusion process.
\newblock {Control $\&$ Cybernetics}  vol. 51, No. 1, 31--42. DOI: 10.2478/candc-2022-0003

\bibitem {abundo:saa20IFPP}
Abundo, M., 2020. \newblock  An inverse problem for the first-passage place of some diffusion processes with random starting point.
\newblock {Stochastic Anal. Appl.}  vol. 38, No. 6, 1122--1133. https://doi.org/10.1080/07362994.2020.1768867

\bibitem {abundo:saa19}
Abundo, M., 2019. \newblock  An inverse first-passage problem revisited: the case of fractional Brownian motion, and time-changed Brownian motion.
\newblock {Stochastic Anal. Appl.}  vol. 37, No. 5, 708--716,  https://doi.org/10.1080/07362994.2019.1608834

\bibitem {abundo:mathematics}
Abundo, M., 2018. \newblock The Randomized First-Hitting Problem
of Continuously Time-Changed Brownian Motion.
\newblock {Mathematics} 6(6), 91, 1--10.
https://doi.org/10.3390/math6060091

\bibitem {abundodelvescovo:mcap17}
Abundo, M. and Del Vescovo, D, 2017. \newblock On the joint distribution of first-passage time and first-passage area of  drifted Brownian motion.
      Methodol Comput Appl Probab  19:985--996 DOI 10.1007/s11009-017-9546-7

\bibitem {abundo:LNSIM}
Abundo, M., 2015. \newblock An overview on inverse first-passage-time problems for one-dimensional diffusion processes.
\newblock {Lecture Notes of Seminario Interdisciplinare di Matematica} Vol. 12, 1 -- 44. http://dimie.unibas.it/site/home/info/documento3012448.html


\bibitem {abundo14b}
Abundo, M., 2014. \newblock One-dimensional reflected diffusions with two boundaries and an inverse first-hitting problem.
\newblock{Stochastic Anal. Appl.} 32, 975-991. DOI: 10.1080/07362994.2014.959595


\bibitem {abundo:MCAP2013}
Abundo, M., 2013. \newblock  On the first-passage area of a one-dimensional jump-diffusion process.
\newblock {Methodol Comput Appl Probab} 15:85--103

\bibitem {abundo13b}
Abundo, M., 2013. \newblock Solving an inverse first-passage-time problem for Wiener process subject to
      random jumps from a boundary.
\newblock{Stochastic Anal. Appl.} 31: 4, 695--707.



\bibitem  {abundo:stapro13}
Abundo, M., 2013. \newblock The double-barrier inverse first-passage problem for Wiener process with random
starting point.
\newblock{Stat. and Probab. Letters} 83, 168--176.

\bibitem  {abundo:stapro12}
Abundo, M., 2012. \newblock An inverse first-passage problem for one-dimensional diffusion with random starting point.
\newblock{Stat. and Probab. Letters} 82, 7--14. \ Erratum: Stat. and Probab. Letters, 82(3), 705.


\bibitem  {abundo:pms00}
Abundo, M., 2000. \newblock On first-passage-times for one-dimensional jump-diffusion processes.
\newblock{Prob. Math.Statis.} 20(2), 399--423.

\bibitem  {borodin:1996}
Borodin, A.N., Salminen, S., 1996. \newblock Handbook of Brownian motion-facts and formulae. Birkhauser Verlag
Basel, Basel

\bibitem  {darling:1953}
Darling, D.A., Siegert, A.J.F., 1953. \newblock The first passage problem for a continuous Markov process.
\newblock {Ann Math Stat}
24:624--639.

\bibitem  {dubey}
Dubey, A. and  Pal, A. (2023) \newblock First-passage functionals for Ornstein Uhlenbeck process with stochastic resetting.
\newblock {J Phys A: Math Theor} 56 (1--19):435002. https://doi.org/10.1088/1751-8121/acf748

\bibitem  {dicre:03}
Di Crescenzo, A., Giorno, V., and  Nobile, A.G. 2003. \newblock On the M/M/1 Queue with Catastrophes and Its Continuous Approximation.
\newblock {Queueing Systems} 43, 329--347.

\bibitem  {gihman:1972}
Gihman, I.I., Skorohod, A.V., 1972. \newblock Stochastic differential equations. Springer, Berlin

\bibitem  {guoyan:24}
Guo, W., Yan, H., and  Chen, H. 2024. \newblock Extremal statistics for for first-passage trajectories of drifted Brownian motion under stochastic resetting.
\newblock {J. Stat. Mech.} 023209, 1--19. https://doi.org/10.1088/1742-5468/ad2678.

\bibitem {gupta}
Gupta, A.K. (Ed.), Nadarajah, S. (Ed.), 2004. \newblock Handbook of Beta Distribution and Its Applications.
\newblock Boca Raton: CRC Press. https://doi.org/10.1201/9781482276596.

\bibitem  {huang:24}
Huang, F.,  and  Chen, H. 2024. \newblock Extremal value statistics of first-passage trajectories of resetting Brownian motion in an interval.
\newblock {J. Stat. Mech.} 093212, 1--21. https://doi.org/10.1088/1742-5468/ad7852.

\bibitem {jackson:stapro09}
Jackson, K., Kreinin, A., and Zhang, W., 2009. \newblock
Randomization in the first hitting problem.
\newblock {Stat. and Probab. Letters} 79, 2422--2428.



\bibitem {karlin2}
Karlin, S., Taylor, H.M., 1981. \newblock A Second Course in Stochastic Processes \newblock
Elsevier.

\bibitem {klebaner}
Klebaner, F.C., 2005. \newblock Introduction to Stochastic Calculus with Applications, 2nd ed. \newblock
London, Imperial College Press.



\bibitem {lanska:89}
Lanska, V. and Smiths C.E., 1989. \newblock The effect of a random initial value in neural first-passage-time
models. \newblock {Math. Biosci.}  93, 191--215.

\bibitem {lefeb:19}
Lefebvre, M., 2019. \newblock Moments of First-Passage Places for Jump-Diffusion Processes.
\newblock {Sankhya A}, 1--9. https://doi.org/10.1007/s13171-019-00181-4

\bibitem {lefeb:22}
Lefebvre, M., 2022. \newblock The inverse first-passage-place problem for Wiener processes.
\newblock {Stochastic Anal. Appl.}, vol. 40, (1), 96--102.

\bibitem{maj07}
Majumdar, S.N., 2007. \newblock Brownian functionals in physics and computer science. \newblock In The Legacy Of
Albert Einstein: A Collection of Essays in Celebration of the Year of Physics (pp. 93-129).

\bibitem {norisa:85}
Nobile, A.G., Ricciardi, L.M., and Sacerdote, L., 1985. \newblock
Exponential trends of Ornstein-Uhlenbeck first-passage-time densities.
\newblock{\it J. Appl. Prob.} {\bf 22}, 360--369.

\bibitem {singh}
Singh, P. and Pal, 2022. \newblock
First-passage Brownian functionals with stochastic
resetting.
\newblock  {J. Phys. A: Math. Theor.} 55 234001: 1--25. https://doi.org/10.1088/1751-8121/ac677c


\bibitem {tuckwell:76}
Tuckwell, H.C., 1976. \newblock On the first-exit time problem for temporally homogeneous Markov processes.
\newblock {J. Appl. Probab.} 13, 39--48.




\end{thebibliography}
\end{document}